\documentclass[preprint,sort&compress,12pt]{elsarticle}
\usepackage{amsmath}
\usepackage{mathrsfs}
\usepackage{stmaryrd}
\usepackage{mathrsfs}
\usepackage{bm}
\usepackage{amsmath}
\usepackage{amsfonts}
\usepackage{amssymb}
\usepackage{amsthm}

\journal{\quad}
\newcommand{\mf}{\mathbf}
\newcommand{\mm}{\mathrm}

\newcommand{\dxdt}{\mm{d} x\mm{d}\tau}
\newcommand{\dydt}{\mm{d}y\mm{d}\tau}

\begin{document}
\begin{frontmatter}
\title{Nonlinear Thermal Instability in Compressible\\ Viscous Flows without Heat Conductivity}
\author[FJ]{Fei Jiang
 }
\ead{jiangfei0591@163.com}
 \address[FJ]{College of Mathematics and Computer Science, Fuzhou University, Fuzhou, 350108, China.}


\begin{abstract}
We investigate the thermal instability of a smooth equilibrium state, in which the density function satisfies Schwarzschild's (instability) condition,
to a compressible heat-conducting viscous flow without heat conductivity in the presence of
a uniform gravitational field in a three-dimensional bounded domain.
We show that the equilibrium state is linearly unstable by a modified variational method. Then, based on the constructed linearly unstable solutions and a local well-posedness
result of classical solutions to the original nonlinear problem, we further construct the initial data of linearly unstable solutions
to be the one of the original nonlinear problem, and establish an appropriate energy estimate of Gronwall-type. With the help of
the established energy estimate, we finally show that the equilibrium state is nonlinearly unstable in the sense of Hadamard
by a careful bootstrap instability argument.
\end{abstract}
\end{frontmatter}

\newtheorem{thm}{Theorem}[section]
\newtheorem{lem}{Lemma}[section]
\newtheorem{pro}{Proposition}[section]
\newtheorem{cor}{Corollary}[section]
\newproof{pf}{Proof}
\newdefinition{rem}{Remark}[section]
\newtheorem{definition}{Definition}[section]

\section{Introduction}
\label{Intro} \numberwithin{equation}{section}

Thermal (or convective) instability often arises when a fluid is heated from below. The classic example of this is a horizontal layer of fluid with its lower side hotter than its upper. The basic state is then one of rests with light and hot fluid below heavy and cool fluid.
When the temperature difference across the layer is great enough, the stabilizing effects of viscosity and thermal conductivity are overcome by the destabilizing buoyancy, and an overturning instability ensues as thermal convection:
hotter part of fluid is lighter and tends to rise as colder part tends to sink according
to the action of the gravity force \cite{DPGRWHHC}.
The phenomenon of thermal convection itself had been recognized by
Rumford \cite{CROFTCW} and  Thomson \cite{TJOACteins}.  
However, the first quantitative experiments on thermal instability and the recognition of the role of viscosity in the phenomenon are due to Benard \cite{BHERESS}, so the convection in a horizontal layer of a fluid heated from below is called B\'enard convection.

The B\'enard convection can be modeled by a (nonlinear) compressible Navier-Stokes-Fourier (simplified by NSF) equations,
in which the coefficients of viscosity and heat conduction are non-zero,
see \cite{GGMPO} for example and the references cited therein.
In 1906, Schwarzschild first derived the criterion for thermal stability of a fluid layer in hydrostatic equilibrium on the basis of Archimedes' principle without considering the effect of viscosity and heat conduction \cite{KSNKGWG,CSAiTTS}. Thus such criterion is named Schwarzschild's criterion.
However, Schwarzschild's argument is at most a suggestive one, since the buoyancy principle applies to the static, not to the perturbed state. Later the criterion was examined rigorously for various situations of the linearized Euler-Fourier model (in other words, the hydrodynamic equations for adiabatic motion, which are linearized for the case of small perturbations from equilibrium)
by other authors, such as Lebowitz \cite{LEbovitzNROth}, Kaniel-Kovetz \cite{KSKA201710110936}, Rosencrans \cite{SIAMJAMOR} and so on.

Unfortunately, the rigorous derivation of Schwarzschild's criterion for the compressible Euler-Fourier model (hyperbolic type) can not be
applied to the compressible NSF model (hyperbolic-parabolic coupled type). Therefore, to investigate the effect of viscosity
and heat conduction in the B\'enard convection, we adopt the
(nonlinear) Boussinesq (approximation) equations, in which the density is considered as a constant in all the terms of the equations
except for the gravity term that is assumed to vary linearly with the temperature \cite{DPGRWHHC}.
It should be noted that this approximation is reasonable, only if the thickness of the layer is small.
Compared with the compressible NSF model, the Boussinesq model enjoys relatively good structure, and moreover, it is parabolic due to
  the absence of the continuity equation of mass conservation (hyperbolic type).
 Thus, based on the linearized Boussinesq model, Rayleigh first showed that instability would occur only when the adverse temperature gradient
 is so large that the Rayleigh number exceeds a certain critical value \cite{RALOC}. Rayleigh's criterion rigorously verifies
the destabilizing effect of buoyancy to the stabilizing effects of viscosity and heat conduction, please refer to the classical monograph \cite{CSHHSCPO} for the details and other results involving the inhibitive effects of the rotation \cite{GGPTA,GGPSBAP} and the magnetic field \cite{GGPPMFR,GGPPMNA,GGPNA} on the thermal convection.
 Moreover, by the energy method and the bootstrap instability method, one can see that Rayleigh's criterion on instability/stability
in the Hadamard sense still holds for the nonlinear Boussinesq model, see \cite{JDDN1966,GYYQH} for example. In addition, the
mathematical theory of attractor  bifurcation  for two-dimensional Boussinesq model \cite{MTWSB} have been established Ma and Wang. However, the corresponding three-dimensional case is still an open problem.

When the density varies by many orders of magnitude (e.g., across a stellar convection zone \cite{MSCTAJ}), the Boussinesq model obviously  fails. Hence, for the B\'enard convection in compressible atmosphere,
 we must adopt the compressible NSF model, and expect to  gain a complete understanding of the effects of density variation and compressibility.
Next we briefly review some results on the B\'enard convection based on the comprehensible NSF model (called compressible B\'enard problem for simplicity).
The results of linear stability and instability  can be founded in \cite{JHTICFHB,SEAAPKC,UWKSMMC,GDOMDRSEAWNO}. Numerical simulations
 for the compressible B\'enard convection are given in \cite{GENJFM} and \cite{CFBNHTJMAHAE} for two- and three-dimensional cases, respectively. Moveover, Gough, Moore, Spiegel and Weiss have noted that the stability or instability depends on the boundary conditions, the heat conduction coefficient and viscosity coefficient from the numerical computations.
 The first theoretical rigorous proof of nonlinear stability was given by Padula and Bollettmo \cite{OMBUMINst}. Successively, Coscia and Padula
further provided a computable critical number and showed a nonlinear stability result for the compressible B\'enard problem
 whenever the Rayleigh number does not exceed the critical number \cite{CosciaPadula}.  Later, Aye and Nishida used the approach of
 a computer assisted proof to obtain the critical Rayleigh number for instability \cite{APTNHCC}.
 Recently, Nishida, Padula and Teramoto proved the existence of steady solutions of the compressible B\'enard problem \cite{TNMPYTHCI,TNMPYTHCII}.

  In \cite{GGMPO}, Guidoboni and Padula
  showed the stability of the linearized compressible B\'enard problem without thermal conductivity under the Schwarzschild's stability  condition.
  At present, it is still open whether the linear result obtained by Guidoboni and Padula can be generalized to the nonlinear case.
  However, in this article,
 we  show that under Schwarzschild's instability condition,
 the nonlinear compressible B\'enard problem without thermal conductivity is unstable.
 We mention that Schwarzschild's instability criterion can be derived on the basis of linearized energy method,
however, our result first rigorously verifies the validity of Schwarzschild's criterion for the nonlinear case.
Next, we shall formulate our problem in details.

\subsection{Rayleigh-B\'enard problem}
The motion of a three-dimensional (3D) compressible
viscous fluid without heat conductivity in the presence of a uniform gravitational field in a
bounded domain $\Omega\subset {\mathbb R}^3$ with smooth boundary is governed by the following equations:
\begin{equation}\label{0101}\left\{\begin{array}{l}
 \rho_t+\mm{div}(\rho{v})=0,\\[1mm]
\rho v_t+\rho {v}\cdot\nabla {v}+\nabla p=\mu\Delta v+\mu_0\nabla\mm{div} v-\rho g{e}_3,\\[1mm]
\rho e_t+\rho v\cdot\nabla e+p\mm{div}v=
{\mu}|\nabla v+\nabla v^\mm{T}|^2/2+\varsigma(\mm{div}v)^2.
\end{array}\right.\end{equation}
Here the unknowns $\rho:=\rho(t, x)$, $v:= v(t, x)$, $e:=e(t, x)$ and $p:=p(x,t)$ denote the
velocity, specific internal energy and pressure of the compressible fluid respectively. $\mu>0$ is the coefficient of shear viscosity,
and $\mu_0=\mu+\varsigma$ with $\varsigma$ being the bulk viscosity, satisfying
$3\varsigma+2\mu\geq 0$. $g>0$ is the gravitational constant, ${e}_3=(0,0,1)^{\mm{T}}$ denotes the vertical unit vector.
The superscript $\mm{T}$ denotes the transposition.

As in \cite{GGMPO}, we assume the pressure satisfies the equations of state for idea gases \cite{FGBFE}, i.e.,
\begin{equation}   \label{201607241416}
p= {R} \rho T_{e}/M,
\end{equation}
where $M$ denotes the molar weight of the gas, $T_{e}$ the (absolute) temperature, and $R$ the ideal gas constant that is equal to the product
of the Boltzmann constant and the Avogadro constant. The internal energy satisfies
\begin{equation}\label{201607241417}        e=c_{V}T_e,
\end{equation}
where ${c_{V}}$ denotes the constant-volume specific heat per mole.��
Under the relations \eqref{201607241416} and \eqref{201607241417}, the equations \eqref{0101} are complete. If we add the dissipative term $\kappa\Delta T_e$
to the right hand side of $\eqref{0101}_3$, then the resulting equations are the Navier-Stokes-Fourier equations with non-zero heat-conduction coefficient $\kappa$.

In this paper we consider the problem of nonlinear convective instability for the equations \eqref{0101} around some equilibrium-state solution,
the density in which satisfies Schwarzschild's instability condition. Thus, we choose an (equilibrium-state) density profile
$\bar{\rho}:=\bar{\rho}(x_3)$, which is independent of $(x_1,x_2)$ and satisfies
\begin{eqnarray}\label{0102}
\bar{\rho}\in C^{4}(\bar{\Omega}),\quad \inf_{ x\in \Omega}\bar{\rho}>0,
\end{eqnarray}
and Schwarzschild's (instability) condition
\begin{eqnarray}\label{01020719}
 -\bar{\rho}'(x_3^0)<\frac{g\bar{\rho}^2(x_3^0)}{(1+a)\bar{p}(x_3^0)}
\mbox{ for some }x^0_{3}\in \{x_3~|~(x_1,x_2,x_3)^{\mm{T}}\in \Omega\},
\end{eqnarray}
where $\bar{\rho}':=\mm{d}\bar{\rho}/\mm{d}x_3$, $a:=R/c_V M$ and $x_{3}^0$ denotes the third component of $x_0\in \Omega$.
We remark that the first condition in \eqref{0102} guarantees that we can construct an unstable classical solution, while
the second one in \eqref{0102} prevents us from treating vacuum in the construction of unstable solutions. Schwarzschild's condition
 results in the classical thermal instability  as will be shown in Theorem \ref{thm:0102} below. To clearly see the physical mechanism, next we
reformulate  \eqref{01020719} by the (equilibrium-state) temperature profile.

In view of the theory of first-order linear ODEs, for given $\bar{\rho}$ in \eqref{0102} we can find a corresponding
(equilibrium-state) internal energy profile $\bar{e}$ that only depends on $x_3$ and is unique up to a constant divided by $\bar{\rho}$, i.e.,
$$\bar{e}=g(a{\bar{\rho}})^{-1}(C- F(\bar{\rho})),$$
where $F(\bar{\rho} )$ denotes some  primitive function of $\bar{\rho} $. In view of \eqref{0102}, we can choose a sufficiently large constant $C$ such that $\bar{e}>0$ on $\bar{\Omega}$.
 Hence, we can construct an internal energy profile $\bar{e}$, such that
\begin{equation}\label{0104}
\bar{e}>0\;\;\mbox{ on }\bar{\Omega},\quad \bar{e}\in C^4(\bar{\Omega}),\;\;\mbox{ and }\;\; \nabla \bar{p}=-\bar{\rho}g  {e}_3\;\; \mbox{ in }\Omega\quad
\mbox{with }\; \bar{p}:=a\bar{\rho}\bar{e}.
\end{equation}
Clearly, the function triple  $(\bar{\rho},0,\bar{e})$ gives an equilibrium state solution to the system \eqref{0101}.
By virtue of the equilibrium state \eqref{0104}, Schwarzschild's condition \eqref{01020719} is equivalent to
\begin{equation}  \label{201607180855}
\bar{e}'(x_3^0)<-\frac{g}{(1+a)}
\mbox{ for some }x^0_{3},
\end{equation}
i.e., the convective condition in view of \eqref{201607241417}:
  \begin{equation}  \label{201607180855new1}
\bar{T}_e'(x_3^0)<-\frac{g}{( c_{V}+R/M)} \mbox{ for some }x^0_{3}.
\end{equation}
From \eqref{201607180855new1} we can easily see why Schwarzschild's condition leads to the convective instability from the physical viewpoint.

We mention that, by Mayer's formula $c_{p}-c_{V}=R$, the convective condition reduces to
\begin{equation}  \label{201607180855new}
\bar{T}_e'(x_3^0)<- {g}/{c_{p} } \mbox{ for some }x^0_{3}
\end{equation}
for $M=1$, where $c_{p}$ denotes the constant-pressure heat capacity per mole. In \cite{GGMPO}, Guidoboni and Padula considered a
layer domain $\mathbb{R}^2\times (0,l)$ with $\Theta_1$ and $\Theta_2$ being the temperatures of the lower and the upper planes respectively, where $\Theta_1>\Theta_2$.
They chose the density profile $\bar{\rho}=\alpha \zeta^m$ and the temperature profile $\bar{T}_e=\beta l\zeta$, where
$$\zeta:=\frac{x_3}{l}+\frac{\Theta_1-\Theta_2}{ \beta l},\ m:=\frac{g}{R\beta l}-1\mbox{ and } \beta:=\frac{\Theta_1-\Theta_2}{d}.$$
Thus, the condition \eqref{201607180855new} further reduces to the condition
$\beta<  {g}/{c_{p} }$, which is a special case of the convective condition \eqref{201607180855new1}.  We mention that
Schwarzschild's condition has another version in terms of entropy, see \cite{LEbovitzNROth,SIAMJAMOR} for example.

Now, we define the perturbation of $(\rho ,v,e)$ around the equilibrium state $(\bar{\rho},0,\bar{e})$ by
$$ \varrho=\rho -\bar{\rho},\quad  u= v- {0},\quad \theta=e-\bar{e}. $$
Then, the triple $(\varrho , u,\theta)$ satisfies the perturbed equations
\begin{equation}\label{0105}\left\{\begin{array}{l}
\varrho_t+\mm{div}((\varrho+\bar{\rho}){ u})=0, \\[1mm]
(\varrho+\bar{\rho}){ u}_t+(\varrho+\bar{\rho}){  u}\cdot\nabla
{  u}+a\nabla [ {({\varrho}+\bar{\rho})(\theta+\bar{e})}-{\bar{\rho}\bar{e}}]=\mu\Delta{ {{u}}}
+\mu_0\nabla\mm{div} u-g\varrho  {e}_3,\\[1mm]
\theta_t+ u\cdot\nabla
(\theta+\bar{e})+a(\theta+\bar{e})\mm{div} u=\{ \mu |\nabla  u+\nabla  u^\mm{T}|^2/2
+\lambda(\mm{div} u)^2\}/(\varrho+\bar{\rho}). \end{array}\right.  \end{equation}
To complete the statement of the perturbed problem, we specify the
initial and boundary conditions:
\begin{eqnarray}\label{0106}   &&
(\varrho,{  u},\theta)|_{t=0}=(\varrho_0,{ u}_0,\theta_0)\quad\mbox{in } \Omega , \\
&& \label{0107}   u|_{\partial\Omega}=0\quad \mbox{ for any }t>0.
\end{eqnarray}
We call the initial-boundary problem \eqref{0105}--\eqref{0107}
the (compressible) Rayleigh-B\'enard problem (without heat conduction). 
For classical solutions of the Rayleigh-B\'enard problem, the initial data should further satisfy the compatibility condition on boundary:
$$ \{(\varrho_0+\bar{\rho}){u_0}\cdot\nabla
{u_0}+a\nabla [ {({\varrho_0}+\bar{\rho})(\theta_0 +\bar{e})}-{\bar{\rho}\bar{e}}]\}|_{\partial\Omega}=(\mu\Delta{u_0}
+\mu_0\nabla\mm{div} u_0 -g\varrho_0 {e}_3)|_{\partial\Omega}. $$
Finally, if we linearize the equations \eqref{0105} around the equilibrium state
$(\bar{\rho}, {0},\bar{e})$, then the resulting linearized equations read as
\begin{equation}\label{0108}
\left\{\begin{array}{ll}
 \varrho_t+\mm{div}( \bar{\rho} u)=0, \\[1mm]
  \bar{\rho} u_t +a\nabla(\bar{e} {\varrho}
  +\bar{\rho}\theta)=\mu\Delta{ {{u}}}+\mu_0\nabla \mm{div} u-g\varrho {e}_3,\\[1mm]
\theta_t+
\bar{e}'u_3+a\bar{e}\mm{div} u = 0.
\end{array}\right.\end{equation} Then the initial-boundary problem \eqref{0106}--\eqref{0108} is called the linearized (compressible) Rayleigh-B\'enard problem (without heat conduction).

\subsection{Main result}
Before stating the main result of this paper, we explain the notations used throughout this paper.
We always assume that the domain $\Omega$ is bounded with smooth boundary.
 For simplicity,
we drop the domain $\Omega$ in Sobolve spaces and the corresponding norms as well as in
integrands over $\Omega$, for example,
\begin{equation*}  \begin{aligned}&
L^p:=L^p(\Omega),\quad
{H}^1_0:=W^{1,2}_0(\Omega),\;\;
{H}^k:=W^{k,2}(\Omega),\;\; \int:=\int_\Omega .  \end{aligned}\end{equation*}
In addition, a product space $(X)^n$ of vector functions
is still denoted by $X$, for examples, the vector function $ u\in (H^2)^3$ is denoted
by $ u\in H^2$ with norm $\| u\|_{H^2}:=(\sum_{k=1}^3\|u_k\|_{H^2}^2)^{1/2}$.

Now, we are able to state our main result on the nonlinear
convective instability in the Rayleigh-B\'enard problem  \eqref{0105}--\eqref{0107}.
\begin{thm}\label{thm:0102}
Under the assumptions \eqref{0102}--\eqref{0104}, the equilibrium state $(\bar{\rho}, 0,\bar{e})$
 is unstable in the Hadamard sense, that is, there are positive constants $\Lambda$, $m_0$, $\varepsilon$ and $\delta_0$,
 and functions $(\bar{\varrho}_0,\bar{ u}_0,\bar{\theta}_0,{ u}_\mm{r})\in H^3$,
such that for any $\delta\in (0,\delta_0)$ and the initial data
 $$ (\varrho_0, u_0,\theta_0):=\delta(\bar{\varrho}_0,\bar{u}_0,\bar{\theta}_0)
 +\delta^2(\bar{\varrho}_0,u_\mm{r},\bar{\theta}_0)\in H^3, $$
there is a unique classical solution $({\varrho},u,\theta)\in C^0([0,T^{\max}),H^3)$ to the Rayleigh-B\'enard problem  \eqref{0105}--\eqref{0107}
satisfying
\begin{equation}\label{0115}
\|(u_1,u_2)(T^\delta)\|_{L^2},\ \|{u}_3(T^\delta)\|_{L^2}\geq {\varepsilon}\;
\end{equation}
for some escape time $T^\delta:=\frac{1}{\Lambda}\mm{ln}\frac{2\varepsilon}{m_0\delta}\in
(0,T^{\max})$, where $T^{\max}$ denotes the maximal time of existence of the solution
$(\varrho, u,\theta)$, and $u_i$ denotes the $i$-th component of $u=(u_1,u_2,u_3)^\mm{T}$.
\end{thm}
\begin{rem}Theorem \ref{thm:0102} still holds for a horizontal periodic domain with finite height,
i.e.,
$$
\Omega:=\{ x:=(x_1,x_2,x_3)^{\mm{T}}\in \mathbb{R}^3~|~(x_1,x_2)\in \mathcal{T},\ 0<x_3<l\}\;\; \mbox{ with }l>0,
$$
  where  $\mathcal{T} :=(2\pi L_1\mathbb{T})\times(2\pi L_2\mathbb{T}) $, $\mathbb{T}=\mathbb{R}/\mathbb{Z}$, and $2\pi L_1$, $2\pi L_2>0$
  are the periodicity lengths. As mentioned before, it is stll an open problem whether the initial-boundary problem \eqref{0105}--\eqref{0107} is stable, provided the density profile satisfies the stability condition
\begin{equation}\label{08172006}
-\bar{\rho}' > {g\bar{\rho}^2 }/{(1+a)\bar{p} }\mbox{ in }\bar{\Omega}.
\end{equation}In a forthcoming ariticle, we will show that the
  corresponding magnetic Rayleigh-B\'enard problem is stable under condition \eqref{08172006}.
\end{rem}
\begin{rem}\label{strongconden}
If the density profile satifies
\begin{equation}\label{dengeq0}    \bar{\rho}'\geq 0,   \end{equation}
which automatically implies Schwarzschild's condition,
then we can establish the instability of the perturbed density, i.e., Theorem \ref{thm:0102} holds with $\|\varrho(T^\delta)\|_{L^2}\geq {\varepsilon}$.
The additional condition \eqref{dengeq0} is used to show $\tilde{ \varrho}_0:=\mm{div}(\bar{\rho}\tilde{v}_0)\not\equiv 0$
in the construction of a linear unstable solution (cf. \eqref{qh0209}), where
$(\tilde{\varrho}_0,\tilde{ v}_0)$ is a solution to the time-independent system \eqref{0109}.��
\end{rem}

Next, we sketch the main idea in the proof of Theorem \ref{thm:0102}. The proof is broken up into three steps.
Firstly,  we make the following ansatz of growing mode solutions to the linearized problem:
\begin{equation}\label{ansatzmode}
(\varrho (x,t),u(x,t),\theta (x,t))=e^{\Lambda t}(\tilde{\rho}(x),\tilde{v}(x),\tilde{\theta}(x)) \quad\mbox{for some }\Lambda>0
\end{equation}
 and deduce \eqref{0108} thus into a time-independent PDE system on the unknown function $\tilde{ v}$.
 Then we adapt and modify the modified variational method in \cite{GYTI2} to the time-independent system to
 get a non-trivial solution $\tilde{ v}$ with a sharp growth rate $\Lambda$, which immediately implies
 that the linearized problem has a unstable solution in the form \eqref{ansatzmode}. This idea was used probably first by
Guo and Tice to deal with an ODE problem arising in constructing unstable solutions of the linearized problem, and later adapted by other researchers
to treat other linear instability problems of viscous fluids, see \cite{JFJSWWWOA,JJTIIA}. Here we directly adapt
this idea to the time-independent PDE system to avoid the use of the Fourier transform and to relax the restriction on domains.
Secondly, we establish the energy estimates of Gronwall-type in $H^3$-norm. Similar (global in time) estimates were obtained
 for the compressible Navier-Stokes-Fourier equations with heat conductivity under the condition of small initial data
 and external forces \cite{MANTIC481,MANTIJ321}. Here we have to modify the arguments in \cite{MANTIC481,MANTIJ321} to deal
 with the compressible Navier-Stokes equations without heat conductivity. Namely, we deal with the sum $\bar{e}\varrho+\bar{\rho}\theta$
as one term (see \eqref{asnoenet}) instead of dividing it into two terms in \cite{MANTIC481}; and we use the equations
\eqref{0105}$_1$ and \eqref{0105}$_2$ independently to control $\|\varrho\|_{H^3}$ and $\|\theta\|_{H^3}$ (i.e. Lemma \ref{lem:0301}),
rather than coupling the equations together to control $\|\varrho\|_{H^3}$ in \cite{MANTIC481}. With these modifications
in techniques, we can get the desired estimates.
Finally, we use the version of the bootstrap instability approach in \cite{GYHCSDDC}
(the interested reader is referred to \cite{FSSWVMNA,FrrVishikM,GYSWIC} for different versions of the bootstrap instability approach)
to show Theorem \ref{thm:0102}, but have to circumvent two additional difficulties due to presence of boundary
which do not appear for spatially periodic problems considered in \cite{GYHCSDDC}:
(i) The idea of Duhamel's principle on the solution operator for linear instability in \cite{GYHCSDDC}
can not be directly applied to our boundary value problem here, since the nonlinear term in \eqref{0105}$_2$ does not vanish
on boundary. To overcome this difficulty, we employ some specific energy estimates to replace Duhamel's principle
(see Lemma \ref{erroestimate} on the error estimate for $\|(\varrho^{\mm{d}},u^{\mm{d}},\theta^{\mm{d}})\|^2_{L^2}$).
(ii) On the boundary the initial data of the linearized Rayleigh-B\'enard problem may not satisfy the compatibility condition
imposed on the initial data of the corresponding nonlinear Rayleigh-B\'enard problem \eqref{0105}--\eqref{0107}. To circumvent
this difficulty, we employ the elliptic theory to construct initial data of the Rayleigh-B\'enard problem that
satisfy the compatibility condition and are close to the initial data of the linearized problem.

The rest of this paper is organized as follows. In Section \ref{sec:02}
we construct unstable solutions of the linearized problem, while in Section \ref{sec:0312} we deduce the nonlinear energy estimates.
Section \ref{sec:04} is devoted to the proof of Theorem \ref{thm:0102}, and finally, in the appendix
we give a proof of the sharp growth rate of solutions to the linearized problem in $H^2$-norm.

\section{Linear instability}\label{sec:02}
In this section, we adapt the modified variational method in \cite{GYTI2} to construct a solution to the
linearized equations \eqref{0108} that has growing $H^3$-norm in time.
We first make a solution ansatz \eqref{ansatzmode} of growing normal mode.
Substituting this ansatz into \eqref{0108}, one obtains the following time-independent system:
\begin{equation}\label{0109}
\left\{
                              \begin{array}{ll}
\Lambda\tilde{\rho}+\mm{div}(\bar{\rho}\tilde{v})=0,\\[1mm]
\Lambda \bar{\rho}\tilde{v}+a\nabla (\bar{e} \tilde{\rho}
  +\bar{\rho}\tilde{\theta})=\mu\Delta\tilde{ v}+\mu_0\nabla \mm{div}\tilde{v}-g
\tilde{\rho}e_3,\\[1mm]
\Lambda   {\tilde{\theta}}+ \bar{e}'\tilde{v}_3
+a\bar{e}\mm{div}\tilde{v}=
0,\\[1mm]
\tilde{ v}|_{\partial\Omega}= 0.
\end{array}
                            \right.
\end{equation}
Eliminating $\tilde{\varrho}$  and $\tilde{\theta}$, one has
 \begin{equation}\label{0113}
\left\{    \begin{array}{ll}
\Lambda^2 \bar{\rho}\tilde{v}+\nabla[
g\bar{\rho}\tilde{{v}}_3-(1+a)\bar{p}\mm{div}\tilde{v}
  ]=\Lambda\mu\Delta\tilde{ {{v}}}+\Lambda\mu_0\nabla \mm{div}\tilde{v}+(
g\bar{\rho}'\tilde{{v}}_3 +g
\bar{\rho}\mm{div}\tilde{v} ){e}_3,\\[1mm]
\tilde{ v}|_{\partial\Omega}=0,
\end{array}           \right.
\end{equation}
where $\tilde{v}$ denotes the third component of $v$.
In view of the basic idea of the modified variational method, we modify the boundary problem \eqref{0113} as follows.
 \begin{equation}\label{nnn0113}
\left\{
                              \begin{array}{ll}
\lambda^2  \bar{\rho}\tilde{v}+\nabla[
g\bar{\rho}\tilde{{v}}_3-(1+a)\bar{p}\mm{div}\tilde{v}
  ]=s\mu\Delta\tilde{ v}+s\mu_0\nabla \mm{div}\tilde{v}+
(g\bar{\rho}'\tilde{{v}}_3+g
\bar{\rho}\mm{div}\tilde{v}) e_3
,\\[1mm]
\tilde{ v}|_{\partial\Omega}= 0,
\end{array}
                            \right.
\end{equation}where $\lambda:=\lambda(s)$ depends on $s$.
We remark that if $s=\Lambda$ is a fixed point of $\lambda(s)$ (i.e., $\lambda(\Lambda)=\Lambda$), then the problem \eqref{nnn0113} becomes \eqref{0113}.

Now, multiplying \eqref{nnn0113}$_1$ by $\tilde{ v}$ and integrating the resulting identity, we get
\begin{equation}\begin{aligned} \label{js1}
\lambda^2 \int \bar{\rho}|\tilde{v}|^2=&\int  \{g\bar{\rho}'\tilde{{v}}_3^2
+[2g\bar{\rho}\tilde{{v}}_3-(1+a)\bar{p}\mm{div}\tilde{v}]\mm{div}\tilde{v}\}\mm{d}  x\\
&-s\int\left(\mu|\nabla
\tilde{ v}|^2+\mu_0 |\mm{div}\tilde{v}|^2\right)\mm{d} x.
\end{aligned}\end{equation}
We define
$$  E_1(\tilde{ v})=\int  \{g\bar{\rho}'\tilde{{v}}_3^2
+[2g\bar{\rho}\tilde{{v}}_3-(1+a)\bar{p}\mm{div}\tilde{v}]\mm{div}\tilde{v}\}
\mm{d} x,$$
and
$$E_2(\tilde{v})=
\int(\mu|\nabla
\tilde{v}|^2+\mu_0|\mm{div}
\tilde{v}|^2)\mm{d} x.$$
Then the standard energy functional for the problem \eqref{nnn0113} is given by
\begin{equation}\label{0204}E(\tilde{v}):=E(\tilde{v},s):=E_1(\tilde{ v})- sE_2(\tilde{v})\end{equation} with an associated admissible set
\begin{equation}\label{0205}
 \mathcal{A}:=\left\{\tilde{ v}\in
H^1_0~\bigg|~J(\tilde{v}):=\int
                            \bar{\rho}\tilde{ v}^2\mm{d} x=1\right\}.
\end{equation}
Recalling (\ref{js1}), we can thus find $\lambda$ by maximizing
\begin{equation}\label{0206}
\lambda^2:=\sup_{\tilde{ v}\in
\mathcal{A}}E(\tilde{ v}).\end{equation}
Obviously, $\sup_{\tilde{ v}\in\mathcal{A}}E(\tilde{ v})<\infty$ for any $s\geq 0$.

Next we show that a maximizer of (\ref{0206}) exists and that the corresponding Euler-Lagrange equations are equivalent to
(\ref{nnn0113}).
\begin{pro}\label{pro:0201}
 Assume that $(\bar{\rho},\bar{e})$ satisfies \eqref{0102} and \eqref{0104},
then for any but fixed $s>0$, the following assertions hold.
 \begin{enumerate}
    \item[(1)] $E({\tilde{ v}})$
achieves its supremum on $\mathcal{A}$.
   \item[(2)] Let $\tilde{ v}_0$ be
a maximizer and $\lambda$ satisfy \eqref{0206},
then  $\tilde{ v}_0\in H^4$  satisfies
the boundary problem (\ref{nnn0113}) and
  \begin{equation}\label{qh0208}
  (\tilde{ v}^0_{1})^2+(\tilde{v}^0_{2})^2\not\equiv 0, \end{equation}
  where $\tilde{v}^0_{i}$ denotes the $i$-th component of $\tilde{v}_0$.
In addition  \begin{equation}\label{qh0209}
  \mm{div}(\bar{\rho}\tilde{ v}_0) \not\equiv 0, \;\;  \mbox{ provided }
\bar{\rho}'\geq 0.
  \end{equation}
 \end{enumerate}
\end{pro}
\begin{pf}
(1) 
Let $\tilde{ v}_n\in \mathcal{A}$ be a maximizing sequence, then
$E(\tilde{ v}_n)$ is bounded from below. This fact together with (\ref{0205}) implies that
$\tilde{ v}_n$ is bounded in $H^1$. So, there exists a
$\tilde{ v}_0\in H^1\cap\mathcal{A}$ and a subsequence (still denoted by $v_n$ for simplicity), such that
$\tilde{ v}_n\rightarrow \tilde{ v}_0$ weakly in $H^1$ and strongly in $L^2$.
 Moreover, by the lower semi-continuity, one has
\begin{equation*}
\begin{aligned}
\sup_{\tilde{ v}\in \mathcal{A}}E(\tilde{ v})
=& \limsup_{n\rightarrow
\infty}E(\tilde{ v}_n)\\
= &\lim_{n\rightarrow \infty}
\int [g\bar{\rho}'(\tilde{v}^n_{3})^2
 +2g\bar{\rho}\tilde{{v}}^n_{3}\mm{div}\tilde{v}_n]\mm{d}
 x
\\
&
-\liminf_{n\rightarrow \infty}
 \int [(1+a)\bar{p}\mm{div}\tilde{v}_n\mm{div}\tilde{v}_n+s\left(\mu|\nabla
\tilde{v}_n|^2+\mu_0 |\mm{div}\tilde{v}_n|^2\right)]\mm{d} x
\\
\leq & E(\tilde{ v}_0)\leq \sup_{\tilde{ v}\in\mathcal{A}}E(\tilde{ v}),
\end{aligned}\end{equation*}
which shows that $E(\tilde{ v})$ achieves its supremum on $\mathcal{A}$.

(2) To show the second assertion, we notice that since $E(\tilde{ v})$ and $J(\tilde{ v})$
 are homogeneous of degree $2$, (\ref{0206}) is equivalent to
  \begin{equation}\label{0227}
  \lambda^2=\sup_{\tilde{ v}\in {H}^1_0}\frac{E(\tilde{ v})}{J(\tilde{ v})}.
\end{equation}
For any $\tau\in \mathbb{R}$ and $ w\in {H}^1_0$, we take
$\tilde{ w}(\tau):=\tilde{ v}_0+\tau w$. Then (\ref{0227}) gives
  \begin{equation*}E(\tilde{ w}(\tau))-\Lambda^2J(\tilde{ w}(\tau))\leq 0.
\end{equation*}
If we set
$I(\tau)=E(\tilde{ w}(\tau))-\Lambda^2J(\tilde{ w}(\tau))$,
then we see that $I(\tau)\in C^1(\mathbb{R})$, $I(\tau)\leq 0$ for all $\tau\in \mathbb{R}$ and
$I(0)=0$. This implies $I'(0)=0$. Hence, a direct computation leads to
  \begin{equation}\label{weakform}\begin{aligned}
&
\int_\Omega \{s\mu\nabla\tilde{v}_0:\nabla w
+[s\mu_0+(1+a)\bar{p}]\mm{div}\tilde{v}_0\mm{div} w\}\mm{d} x\\
&=\int_\Omega [g\bar{\rho}{\mathrm{div}} \tilde{ v}_0 e_3 + g
\bar{\rho}'\tilde{{v}}^0_{3} e_3 -\nabla(g\bar{\rho}\tilde{v}^0_{3})-\Lambda^2\bar{\rho}\tilde{ v}_0]\cdot\tilde{ w}
 \mm{d} x.\end{aligned}  \end{equation}
which shows that $\tilde{ v}$ is a weak solution to the boundary problem \eqref{nnn0113}.
Recalling that $0<\bar{p}\in C^4(\bar{\Omega})$, $\bar{\rho}\in C^4(\bar{\Omega})$ and $\tilde{v}_0\in H^1(\Omega)$, by a bootstrap argument
and the classical elliptic theory, we infer from the weak form \eqref{weakform} that $\tilde{v}_0\in H^4(\Omega)$.

Next we turn to the proof of \eqref{qh0208} and \eqref{qh0209} by contradiction. Suppose
that $(\tilde{ v}^0_{1})^2+(\tilde{v}^0_{2})^2\equiv 0$ or $\mm{div}(\bar{\varrho}\tilde{v}_{0})\equiv 0$, then
\begin{equation}\label{horizve}\begin{aligned}0<\lambda^2 =&\int  \{g\bar{\rho}'(\tilde{{v}}^0_{3})^2
+[2g\bar{\rho}\tilde{{v}}^0_{3}-(1+a)\bar{p}\partial_{x_3}\tilde{v}^0_{3}]
\partial_{x_3}\tilde{v}_{3}^0\}\mm{d}  x-s\int\left(\mu|\nabla
\tilde{ v}^0_{3}|^2+\mu_0  |\partial_{x_3}\tilde{v}_{3}^0|^2\right)\mm{d} x\\
=&-\int  (1+a)\bar{p}|\partial_{x_3}\tilde{v}_{3}^0|^2\mm{d}  x-s\int\left(\mu|\nabla
\tilde{ v}_{3}^0|^2+\mu_0 |\partial_{x_3}\tilde{v}_{3}^0|^2\right)\mm{d} x< 0,
\end{aligned}
\end{equation}
or
\begin{equation}\begin{aligned} \label{densityproe}
0<\lambda^2 =&\int  \{g\bar{\rho}'(\tilde{{v}}_{3}^0)^2
+[2g\bar{\rho}\tilde{{v}}_{3}^0-(1+a)\bar{p}\mm{div}\tilde{v}_{0}]\mm{div}\tilde{v}_0\}\mm{d}  x-s\int\left(\mu|\nabla
\tilde{ v}_0|^2+\mu_0 |\mm{div}\tilde{v}_0|^2\right)\mm{d} x\\
=&-\int[g\bar{\rho}'(\tilde{{v}}^0_{3})^2
+(1+a)\bar{p}|\mm{div}\tilde{v}_0|^2]\mm{d}  x-s\int\left(\mu|\nabla
\tilde{ v}_0|^2+\mu_0 |\mm{div}\tilde{v}_0|^2\right)\mm{d} x<0 ,
\end{aligned}\end{equation}
which contradicts. Therefore, \eqref{qh0208} and \eqref{qh0209} hold.
 This completes the proof.
\hfill $\Box$
\end{pf}

Next, we want to show that there is a fixed point $\Lambda$ such that $\lambda(\Lambda)=\Lambda>0$. To this end,
we first give some properties of $\alpha(s):=\sup_{\tilde{ v}\in
\mathcal{A}}E(\tilde{ v},s)$ as a function of $s> 0$.
\begin{pro}\label{pro:0202}
 Assume that $(\bar{\rho},\bar{e})$ satisfies \eqref{0102}--\eqref{0104}. Then the
function $\alpha(s)$ defined on $(0,\infty)$ enjoys the following properties:
\begin{enumerate}[\quad \ (1)]
 \item $\alpha(s)\in C_{\mathrm{loc}}^{0,1}(0,\infty)$ is nonincreasing.
   \item
 There are constants $c_1$, $c_2>0$ which depend on $g$, $\bar{\rho}$ and $\mu$, such that
  \begin{equation}\label{0210}\alpha(s)\geq c_1-sc_2
  .\end{equation}
  \end{enumerate}
\end{pro}
\begin{pf}
(1) Let $\{\tilde{ v}^n_{s_i}\}_{n=1}^{+\infty}\subset\mathcal{A}$ be a maximizing sequence
of $\sup_{\tilde{ v}\in\mathcal{A}}E(\tilde{ v},s_i)=\alpha(s_i)$ for $i=1$ and $2$. Then
\begin{equation*}
\alpha(s_1)\geq
\limsup_{n\rightarrow\infty}E(\tilde{ v}_{s_2}^n,s_1) \geq
\liminf_{n\rightarrow\infty}E(\tilde{ v}_{s_2}^n,s_2)=\alpha(s_2)\;
\mbox{ for any }0<s_1<s_2<\infty.
\end{equation*}
Hence $\alpha(s)$ is nonincreasing on $(0,\infty)$. Next we use this fact to show the continuity of $\alpha(s)$.

Let $I:=[b,c]\subset (0,\infty)$ be a bounded interval.
Noting that, by Cauchy-Schwarz's inequality,
\begin{equation*} \label{}
\begin{aligned}
E(\tilde{v})
\leq &\int (g\bar{\rho}'\tilde{{v}}_3^2
+2g\bar{\rho}\tilde{{v}}_3\mm{div}\tilde{v})\mm{d}  x- (1+a)\int\bar{p}|\mm{div}\tilde{v}|^2\mm{d} x  \\
\leq &{g}\left[\left\|\frac{\bar{\rho}'}{\bar{\rho}}
\right\|_{L^\infty}+\frac{g}{(1+a)}\left\|\frac{\bar{\rho}}{\bar{p}}\right\|_{L^\infty}\right].
\end{aligned}\end{equation*}
Hence, by the monotonicity of $\alpha (s)$ we have
\begin{equation} \label{02321}
|\alpha(s)|\leq \max\left\{|\alpha(b)|,{g}\left[\left\|\frac{\bar{\rho}'}{\bar{\rho}}
\right\|_{L^\infty}+\frac{g}{(1+a)}\left\|\frac{\bar{\rho}}{\bar{p}}\right\|_{L^\infty}\right]\right\}:=
L<\infty\quad\mbox{ for any }s\in I.
\end{equation}
On the other hand, for any $s\in I$, there exists a maximizing sequence
$\{\tilde{ v}^n_{s}\}\subset\mathcal{A}$ of $\sup_{\tilde{ v}\in
\mathcal{A}}E(\tilde{ v},s)$, such that
\begin{equation}\label{0232}\begin{aligned}|\alpha(s)-E(\tilde{ v}_{s}^n,s)|<1
\end{aligned}.\end{equation}
Making use of (\ref{0204}), (\ref{02321}) and (\ref{0232}), we infer that
\begin{equation*}\begin{aligned}\label{0234}0\leq&
\int\left(\mu |\nabla \tilde{ v}|^2+\mu_0 |\mm{div} \tilde{ v}|^2\right)\mm{d} x
\\=&\frac{1}{s}\int \{g\bar{\rho}'|\tilde{{v}}_{s3}^n|^2
+[2g\bar{\rho}\tilde{{v}}_{s3}^n
-(1+a)\bar{p}\mm{div}\tilde{v}_s^n]\mm{div}\tilde{v}_s^n\}\mm{d}
 x -\frac{E(\tilde{ v}_s^n,s)}{s}
\\
\leq &
\frac{1+L}{b}+\frac{g}{b}\left[\left\|\frac{\bar{\rho}'}{\bar{\rho}}
\right\|_{L^\infty}+\frac{g}{(1+a)}\left\|\frac{\bar{\rho}}{\bar{p}}\right\|_{L^\infty}\right]
:=K.
\end{aligned}\end{equation*}
Thus, for $s_i\in I$ ($i=1,2$), we further find that
\begin{equation}\begin{aligned}\label{0235}\alpha(s_1)= \limsup_{n\rightarrow
\infty}E(\tilde{ v}_{s_1}^n,s_1)\leq & \limsup_{n\rightarrow
\infty}E(\tilde{ v}_{s_1}^n,s_2)\\
&+ |s_1-s_2|\limsup_{n\rightarrow
\infty}\int_{\Omega}(\mu|\nabla \tilde{ v}_{s_1}^n|^2+\mu_0|\mm{div} \tilde{ v}_{s_1}^n|^2)\mm{d} x\\
\leq & \alpha(s_2)+K|s_1-s_2|.
\end{aligned}\end{equation}

Reversing the role of the indices $1$ and $2$ in the derivation of the inequality
(\ref{0235}), we obtain the same boundedness with the indices switched. Therefore, we deduce that
\begin{equation*}\begin{aligned}|\alpha(s_1)-\alpha(s_2)|\leq K|s_1-s_2|,
\end{aligned}\end{equation*}
which yields $\alpha(s)\in C_{\mathrm{loc}}^{0,1}(0,\infty)$.

(2) We turn to prove \eqref{0210}. First we should construct a function ${v}\in H_0^1$, such that
\begin{equation}\label{0214}
E(v)>0.
\end{equation}
In view of Schvarzschild's  condition,
there is a cylinder $$\Xi_{x_0}^{ \delta}:=\left\{ x\in \mathbb{R}^3 ~\bigg|~\sqrt{( x_1- x_1^0)^2+( x_2- x_2^0)^2}\leq \delta,\ |x_3-x_3^0|\leq \delta\right\}\subset \Omega,$$
such that
$$ -\bar{\rho}' <\frac{g\bar{\rho}^2 }{(1+a)\bar{p} } \mbox{ on }\Xi_{ x_0}^{r,\delta}.$$

Now, choose a smooth function $f(z) \in H^1_0(\mathbb{R})$, such that
\begin{equation*}f(z)\left\{
 \begin{array}{ll}
>0 , & \hbox{ for }  |z| < \delta/2; \\
=0, & \hbox{ for } |z| \geq \delta/2
                   \end{array}
                 \right.
 \end{equation*}
Then we  define
\begin{equation*}\bar{v}(x):=
\left(0,\varphi(x_3)f\left(\sqrt{x_1^2+x_2^2}\right),f(x_3)\partial_2
f\left(\sqrt{x_1^2+x_2^2}\right)\right),
\end{equation*}
where $$\varphi=\frac{g\bar{\rho}}{(1+a)\bar{p}}f-f'.$$
It is easy to check that
${ { v}}:=\bar{ v}( x- x_0)\in H_0^1(\Omega)
$ satisfies
$$ \int  g\left[\bar{\rho}'+\frac{g\bar{\rho}^2}{(1+a)\bar{p}}
\right] {{v}}_3^2\mm{d}  x>0
\mbox{ and }  \frac{g\bar{\rho} {{v}}_3}{ {(1+a)\bar{p}}}
-\mm{div}\tilde{v} =0,$$
which, together with the relation
$$E (v) =\int  \left\{g\left[\bar{\rho}'+\frac{g\bar{\rho}^2}{(1+a)\bar{p}}
\right] {{v}}_3^2
- {(1+a)\bar{p}}\left[\frac{g\bar{\rho} {{v}}_3}{ {(1+a)\bar{p}}}
-\mm{div} {v}\right]^2\right\}\mm{d}x ,$$
 implies that  \eqref{0214}.

With \eqref{0214} to hand, one has
\begin{equation*}\begin{aligned}
\alpha(s)=& \sup_{\tilde{ v}\in
\mathcal{A}}E(\tilde{ v},s)=\sup_{\tilde{ v}\in
{H}^1_0}\frac{E(\tilde{ v},s)}{J(\tilde{ v})}  \\
&\geq \frac{E({  { v}},s)}{J({ { v}})}= \frac{
E_1
(v)}{\int\bar{\rho}{ v}^2\mm{d} x}
-s\frac{ \mu \int|\nabla { v}|^2\mm{d} x}{\int\bar{\rho}{ v}^2\mm{d} x}:= c_1-sc_2
\end{aligned}\end{equation*}
for two positive constants $c_1:=c_1(g,\bar{\rho})$ and $c_2:=c_2(g,\mu,\bar{\rho})$.
This completes the proof of Proposition \ref{pro:0202}.
 \hfill $\Box$
\end{pf}

Next we show that there exists a function $\tilde{{v}}$ satisfying
\eqref{0113} with a grow rate $\lambda$. Let
\begin{equation*}\label{}\mathfrak{S} :=\sup\{s~|~\alpha(\tau)>0\;\mbox{ for any }\tau\in (0,s)\}.
\end{equation*}
By virtue of Proposition \ref{pro:0202}, $\mathfrak{S}>0$; and moreover, $\alpha(s)>0$ for any $s<\mathfrak{S}$.
Since $\alpha(s)=\sup_{\tilde{ v}\in\mathcal{A}}E(\tilde{ v},s)<\infty$, we make use of
the monotonicity of $\alpha(s)$ to deduce that
 \begin{equation}\label{zero}
 \lim_{s\rightarrow 0}\alpha(s)\mbox{ exists and the limit is a positve constant.}
 \end{equation}

 On the other hand, by virtue of Poinc\'{a}re's inequality, there is a constant $c_3$ dependent of
$g$, $\bar{\rho}$ and $\Omega$, such that
 $$
 \begin{aligned}
 g\int( \bar{\rho}'\tilde{v}_3^2+2\bar{\rho}\tilde{{v}}_3\mm{div}\tilde{v})\mm{d} x\leq
c_3
\int|\nabla\tilde{v}|^2\mm{d} x\quad\mbox{ for any }\tilde{ v}\in\mathcal{A}.
\end{aligned}$$
 Thus, if $s>c_3/\mu$, then
 $$g\int(\bar{\rho}'\tilde{{v}}_3^2
+2\bar{\rho}\tilde{{v}}_3\mm{div}\tilde{v})\mm{d}
 x-s\mu\int|\nabla
\tilde{ v}|^2\mm{d} x<0\quad\mbox{ for any }\tilde{ v}\in\mathcal{A},$$
which implies that
 $$\alpha(s)\leq 0\quad \mbox{ for any } s>c_3/\mu. $$
Hence $\mathfrak{S}<\infty$, and moreover,
\begin{equation}\label{zerolin}  \lim_{s\rightarrow \mathfrak{S}}\alpha(s)=0.  \end{equation}

Now, employing a fixed-point argument, exploiting \eqref{zero}, \eqref{zerolin}, and the continuity of $\alpha(s)$ on
$(0,\mathfrak{S})$, we find that there exists a unique $\Lambda\in(0,\mathfrak{S})$, such that
\begin{equation} \label{growth}
 \Lambda=\sqrt{\alpha(\Lambda)}=\sqrt{\sup_{\tilde{ w}\in
\mathcal{A}}E(\tilde{ w}, \Lambda)}>0. \end{equation}

In view of Proposition \ref{pro:0201}, there is a solution $\tilde{v}\in H^4$
to the boundary problem (\ref{nnn0113}) with $\Lambda$ constructed in \eqref{growth}. Moreover,
$\Lambda^2=E(\tilde{v},\Lambda)$, $\tilde{ v}_{1}^2+\tilde{v}_{2}^2\not\equiv 0$ and
$\tilde{{v}}_3\not\equiv 0$ by \eqref{growth} and \eqref{0204}. In addition, $\mm{div}(\bar{\rho}\tilde{v})\not\equiv 0$  provided
$\bar{\rho}'\geq 0$. Thus we have proved
\begin{pro}\label{pro:nnn0203}
 Assume that $(\bar{\rho},\bar{e})$ satisfies \eqref{0102}--\eqref{0104}. Then there exists a
$\tilde{v}\in H^{4}$ satisfying the boundary problem \eqref{0113}
with a growth rate $\Lambda>0$ defined by
\begin{equation}   \label{sharprate}
\Lambda^2=\sup_{\tilde{ w}\in {H}_0^1(\Omega)}\frac{E_1(\tilde{ w})
-\Lambda E_2(\tilde{ w})}{\int\bar{\rho}|\tilde{ w}|^2\mathrm{d} x}.
\end{equation}
 Moreover, $\tilde{v}$ satisfies $\mm{div}(\bar{\rho}\tilde{v})\not\equiv 0$, $\tilde{ v}_{1}^2+\tilde{v}_{2}^2\not\equiv 0$ and $\tilde{{v}}_3\not\equiv 0$.
 In particular, let $(\tilde{\rho},\tilde{\theta}):=-(\mm{div}(\bar{\rho}\tilde{v}), \bar{e}'\tilde{v}_3
+a\bar{e}\mm{div}\tilde{v})/\Lambda$, then
 $(\tilde{\rho},\tilde{v},\tilde{\theta})\in H^3$ satisfies \eqref{0109}.
In addition, $\tilde{\rho}\not\equiv 0$ provided $\bar{\rho}'\geq 0$.
\end{pro}

As a result of Proposition \ref{pro:nnn0203}, one immediately gets the following linear instability.
\begin{thm}\label{thm:0101}
 Assume that $(\bar{\rho},\bar{e})$ satisfies \eqref{0102}--\eqref{0104}.
Then the equilibrium state $(\bar{\rho}, {0},\bar{e})$
 is linearly unstable. That is, there exists an unstable solution
$$(\mf{\varrho}, u,\theta):=e^{\Lambda t}(\tilde{\rho},\tilde{v},\tilde{\theta})$$
 to  the linearized   Rayleigh-B\'enard problem \eqref{0106}--\eqref{0108}, such that $(\tilde{\rho},\tilde{v},\tilde{\theta})\in H^3$ and
\begin{equation*}
\|({u}_1,u_2)(t)\|_{L^2}\mbox{ and } \|{u}_3(t)\|_{L^2}\to \infty\mbox{ as }t\to\infty ,
\end{equation*}
where the constant growth rate $\Lambda$ and $(\tilde{\varrho},\tilde{v},\tilde{\theta})$ are constructed in Proposition \ref{pro:nnn0203}.
Moreover, $\tilde{\rho}\not\equiv 0$ provided $\bar{\rho}'\geq 0$.
\end{thm}

\section{Nonlinear energy estimates}\label{sec:0312}
In this section, we derive some nonlinear energy estimates for the  (nonlinear) Rayleigh-B\'enard problem
\eqref{0105}--\eqref{0107} and an estimate of Gronwall-type in $H^3$-norm, which will  be used in the
proof of Theorem \ref{thm:0102} in the next section. To this end, let $(\varrho,{u},\theta)$ be a solution
of the Rayleigh-B\'enard problem, such that
\begin{equation}\label{enerdienf}
\mathcal{E}(t):=\mathcal{E}(\varrho,u,\theta)(t):=\|(\varrho, u,\theta)(t)\|_{H^3} \leq \delta_1^0,
\end{equation}
where $\delta_0^1$ is sufficiently small. It should be noted that the smallness depends on the physical parameters in \eqref{0105}, and satisfies the following property by
using the embedding $H^3\hookrightarrow L^\infty$:
\begin{equation*}\label{}0<\frac{\inf_{x\in \Omega}\{\bar{\rho}\}}{2}
\leq \rho(t, x):=\varrho+\bar{\rho}\leq 2{\sup_{x\in \Omega}\{\bar{\rho}\}}\mbox{ for any }t\geq 0,\ x\in \Omega,
\end{equation*}
where $\underline{\rho}$ and $\bar{\rho}$ are constants. We remark here that these assumptions will be repeatedly used in what follows.
Moreover, we assume that the solution $(\varrho,{u},\theta)$ possesses proper regularity, so that the procedure
of formal calculations makes sense. For simplicity, we only sketch the outline and shall omit the detailed calculations.
We remind that in the calculations that follow, we shall repeatedly use the Sobolev embedding theorem \cite[Subsection 1.3.5.8]{NASII04},
Young's, H\"older's and Poincar\'e's inequalities, and the following interpolation inequality \cite[Chapter 5]{ARAJJFF}:
 $$\|f\|_{H^j}\lesssim \|f\|_{L^2}^{\frac{1}{j+1}}\|f\|_{H^{j+1}}^{\frac{j}{j+1}}
\leq C_\epsilon\|f\|_{L^2} +\epsilon\|f\|_{H^{j+1}}\mbox{ for any constant }\epsilon>0. $$

In addition, we shall always use the following abbreviations in what follows.
\begin{align}
&\mathcal{E}_0:={\mathcal{E}}(\varrho_0,{u}_0,{N}_0),\  D^k:=\{\partial_{x_{1}}^{k_1}\partial_{x_{2}}^{k_2}
\partial_{x_{3}}^{k_3}\}_{k_1+k_2+k_3=k},\nonumber\\
&|\|g D^k f|\|^2:=\sum_{k_1+k_2+k_3=k}|\|g \partial_{x_{1}}^{k_1}\partial_{x_{2}}^{k_2}\partial_{x_{3}}^{k_3}f|\|^2
\;\mbox{ for some norm }|\|\cdot|\|,\nonumber\\
&\frac{d}{dt}:=\partial_t+u\cdot\nabla \;\mbox{ denotes the material derivative},\nonumber\\
& L^\varrho\equiv L^\varrho(\varrho, u):= \varrho_t+\bar{\rho}'{u}_3
+  \bar{\rho}\mm{div} u,\ N^\varrho:=N^\mf{\varrho}(\varrho, u):=-\mm{div}(\varrho u), \nonumber
\\[1mm]
&
 {L}^u\equiv {L}^u(\varrho, u,\theta):= \bar{\rho} u_t +a\nabla(\bar{e}
 {\varrho}  +\bar{\rho}\theta )-\mu \Delta  u-\mu_0\nabla \mm{div} u+g\varrho  e_3,\nonumber  \\
& {N}^u\equiv {N}^u(\varrho,\theta, u):= -( \varrho+\bar{\rho}) u\cdot\nabla u- \varrho u_t -a\nabla(\varrho \theta),  \nonumber\\[1mm]
&
L^\theta\equiv L^\theta(\varrho, u,\theta):=\theta_t +\bar{e}'{u}_3+a\bar{e}\mm{div} u,\nonumber
\\[1mm]
&N^\theta\equiv N^\theta(\varrho,u,\theta):=[{\mu}|\nabla  u+\nabla ( u)^\mm{T}|^2/2+ \lambda(\mm{div} u)^2]/(\varrho+\bar{\rho})
- u\cdot\nabla\theta-a\theta\mm{div} u,  \nonumber\\
& \mathcal{R}(t):= \left\|\left(\varrho,\theta,u_t,\frac{{d}}{{d}t}\left(\bar{e}\varrho+\bar{\rho} \theta\right)\right)
\right\|_{H^2}^2+\mathcal{E}(\|u\|_{H^3}+\|u\|_{H^4}^2+\mathcal{E}^2),\nonumber\\
 &a\lesssim b\mbox{ means that }a\leq Cb\mbox{ for some constant }C>0,\nonumber\end{align}
where the constant $C$  may depend on some physical parameters in the perturbed equations \eqref{0105}.
In particular, the perturbed equations can be written as the following non-homogenous form:
\begin{align}
&\label{masssimply} L^\varrho=N^\varrho,
   \\
&\label{mometursim}{L}^u =N^u,
\\[1mm]
& \label{masssimply2}L^\theta=N^\theta.\end{align}
In addition, we can use \eqref{masssimply} and \eqref{masssimply2} to
deduce that
\begin{equation}
\label{relationofdiv}\mm{div}u=-\frac{1}{\bar{\rho}\bar{e}+
\bar{p}}\frac{d}{dt}(\bar{e}\varrho+\bar{\rho}\theta)+\frac{(\bar{e}'\varrho +  \bar{\rho}'\theta )u_3+\bar{\rho}(N^\theta +
u\cdot \nabla \theta) -(\bar{\rho} \bar{e})'u_3
- \bar{e}\varrho\mm{div}u}{\bar{\rho}\bar{e}+\bar{p}}.
\end{equation}
Thus
 \eqref{mometursim}   can be rewritten as follows, which will be used in the boundary estimates.
\begin{equation}
 \label{newform}{L}^u_{\mm{new}} =N^u_{\mm{new}},\end{equation}
where
$${L}^u_{\mm{new}}:=  \bar{\rho} u_t -\mu \Delta  u+  g\varrho  e_3
+\nabla \left[\frac{\mu_0}{\bar{\rho}\tilde{e}+\bar{p}}
\frac{d}{dt}(\bar{e} {\varrho}+\bar{\rho} {\theta})+a\bar{e} {\varrho}+a\bar{\rho}
 {\theta}\right]  $$and
$$N^u_{\mm{new}}:=N^u+\mu_0\nabla \{[(\bar{e}'\varrho +  \bar{\rho}'\theta )u_3+\bar{\rho}(N^\theta +
u\cdot \nabla \theta) -(\bar{\rho} \bar{e})'u_3
- \bar{e}\varrho\mm{div}u ]/
 (\bar{\rho}\bar{e}+\bar{p})\}. $$

Next, we shall establish a series of lemmas which imply \emph{a priori } estimates for the perturbed density, velocity and temperature.
\subsection{Estimates on the whole domain}
Firstly, we have the following estimate on the perturbed density and temperature.
\begin{lem}\label{lem:0301}
For $0\leq k\leq 3$, it holds that
$$\| (\varrho,\theta)(t)\|_{H^k}^2\lesssim \|(\varrho,\theta)(0)\|_{H^k}^2+\int_0^t\mathcal{E}(\|u\|_{H^{k+1}}
+\mathcal{E}^2)\mm{d}\tau.$$
\end{lem}
\begin{pf}Using the identity
\begin{equation*}  \int_0^t\int D^k L^{\theta}D^k\theta\dxdt=\int_0^t\int D^k N^\theta D^k \theta\dxdt\quad\mbox{ for } 0\leq k\leq 3,
\end{equation*}
we have
\begin{align*}
   \|D^k\theta(t)\|_{L^2}=&\|D^k\theta(0)\|_{L^2}
   -\int_0^t\int D^k(\bar{e}'{u}_3+a\bar{e}\mm{div} u)D^k\theta\dxdt\nonumber\\
   &-\int_0^t\int D^k(u\cdot\nabla\theta+a\theta\mm{div} u)D^k\theta\dxdt\nonumber\\
   &+\int_0^t\int D^k\{[{\mu}|\nabla  u+\nabla ( u)^\mm{T}|^2/2+ \lambda(\mm{div} u)^2]/(\varrho+\bar{\rho})\}D^k\theta\dxdt
\end{align*}
Noting that $\mathcal{E}\leq 1$ and $$\int u\cdot\nabla  D^k    \theta D^k\theta\mm{d}x=
-\frac{1}{2}\int   | D^k\theta |^2 \mm{div}u\mm{d}x,$$
hence it's easy to get
\begin{align}  \label{newJJ0305}
   \|D^k\theta(t)\|_{L^2}
\lesssim \|D^k\theta(0)\|_{L^2}^2+\int_0^t \|D^k\theta\|_{L^2}(\|u\|_{H^{k+1}}
+\mathcal{E}^2)\mm{d}\tau.
\end{align}

Similarly, using the identity
\begin{equation*}  \int_0^t\int D^k L^{\varrho}D^k\varrho\dxdt=\int_0^t\int D^k N^\varrho D^k \varrho\dxdt \quad\mbox{ for } 0\leq k\leq 3,
\end{equation*}
we arrive at
\begin{equation}\label{newJJ0306}  \| \varrho(t)\|_{H^k}^2\lesssim \|D^k\varrho(0)\|_{L^2}^2+\int_0^t\|D^k\varrho\|_{L^2}(\|u\|_{H^{k+1}}
+\mathcal{E}^2)\mm{d}\tau.
\end{equation}
Summing up \eqref{newJJ0305} and \eqref{newJJ0306}, we immediately get the desired conclusion.
\hfill$\Box$
\end{pf}

Secondly, we control the perturbed velocity.
\begin{lem}\label{lem:0303}   It holds that
\begin{equation*}\label{}\begin{aligned}\|u\|_{H^{3}}^2
\lesssim \| u_t\|_{H^{1}}^2+
 \|(\varrho,\theta)\|_{H^2}^2+\mathcal{E}^4.
\end{aligned}\end{equation*}
\end{lem}
\begin{pf}
Since the viscosity term in \eqref{mometursim} defines a strongly elliptic operator on $u$, we have for
$u\in H^k\cap H_0^1$ ($1\leq k\leq 3$) that
\begin{equation}\label{ellioper}
\|u\|_{H^k}^2\lesssim \|\mu \Delta  u+\mu_0\nabla \mm{div} u\|_{H^{k-2}}^2.
\end{equation}
Thus, applying  \eqref{ellioper} to the system
 \begin{equation}\label{ellipicequation}
-\mu \Delta  u-\mu_0\nabla \mm{div} u ={N}^u-\bar{\rho} u_t-g\varrho  e_3-a\nabla(\bar{e} {\varrho} +\bar{\rho}\theta ),
\end{equation}
one immediately get the desired conclusion.
\hfill$\Box$
\end{pf}

Thirdly, we bound the time-derivative of the perturbed velocity.
\begin{lem}\label{lem:0304} It holds that
\begin{eqnarray}  && \label{detetim1}
\|(\varrho,\theta)_t\|_{H^k}^2\lesssim\| u\|_{H^{k+1}}^2+\mathcal{E}^4
\lesssim \mathcal{E}^2\quad \mbox{ for }0\leq k\leq 2,  \\
&& \label{momentum2}  \|  u_t(t)\|_{H^1}^2+\int_0^t\|u_{\tau\tau}\|^2_{L^2}
\mm{d}\tau \lesssim \|D^1u_t|_{t=0}\|_{L^{2}}^2+ \int_0^t( \|u\|_{H^{2}}^2 +\mathcal{E}^4)\mm{d}\tau , \\
 && \label{utt}  \| u_t\|_{H^{2}}^2\lesssim \| u_{tt}\|_{L^2}^2+ \|u\|_{H^2}^2+\mathcal{E}^4.
 \end{eqnarray}
 \end{lem}
\begin{pf}
The inequality (\ref{detetim1}) follows directly from \eqref{masssimply} and \eqref{masssimply2}.
By \eqref{mometursim} we see that
\begin{equation}\label{detetim2}\|u_t\|_{H^1}^2\lesssim
\|(\varrho,\theta)\|_{H^2}^2+\| u\|_{H^{3}}^2+\mathcal{E}^4\lesssim \mathcal{E}^2.\end{equation}
On the other hand, noting that
$$\int_0^t\int {L}^{ u}_\tau\cdot u_{\tau\tau}\dxdt =\int_0^t\int {N}^{ u}_\tau\cdot u_{\tau\tau}\dxdt,  $$
we  have $$\begin{aligned}
&\mu\|\nabla  u_{t}(t)\|_{L^2}^2+\mu_0\|\mm{div}  u_{t}(t)\|_{L^2}^2+\int_0^t
\|\sqrt{\rho}u_{\tau\tau}\|_{L^2}^2\mm{d}\tau\\
&=\mu\|\nabla  u_{t}|_{t=0}\|_{L^2}^2+\mu_0\|\mm{div}  u_{t}|_{t=0}\|_{L^2}^2-\int_0^t\int [a\nabla(\bar{e}
 {\varrho}  +\bar{\rho}\theta )+g\varrho  e_3]_\tau\cdot u_{\tau\tau}\dxdt\\
&\quad-\int_0^t\int[ ( \varrho+\bar{\rho}) u\cdot\nabla u+ \varrho u_\tau +a\nabla(\varrho \theta)]_\tau\cdot u_{\tau\tau}\dxdt\\
&\leq C \|D^1u_t|_{t=0}\|_{L^2}+C \int_0^t[\|(\varrho, \theta)_\tau\|_{H^1}^2+
 \|(\varrho_\tau,u_\tau,\theta_\tau)\|_{H^1}^2 (\mathcal{E}^2+\|\varrho_\tau \|_{H^1}^2)]\mm{d}\tau\\
&\qquad +\frac{1}{2}\int_0^t
\|\sqrt{\rho}u_{\tau\tau}\|_{L^2}^2\mm{d}\tau +\left\|\frac{ {\varrho}}{\rho}\right\|_{L^\infty}\int_0^t
\|\sqrt{\rho }u_{\tau\tau}\|_{L^2}^2\mm{d}\tau.
\end{aligned}$$Using the embedding $H^2\hookrightarrow L^\infty$ for $\varrho$, and taking $\delta^0_1$ to be sufficiently small, we get
 $$\begin{aligned}
&\mu\|\nabla  u_{t}(t)\|_{L^2}^2+\mu_0\|\mm{div}  u_{t}(t)\|_{L^2}^2+\int_0^t
\|\sqrt{\rho}u_{\tau\tau}\|_{L^2}^2\mm{d}\tau \\
&\lesssim \|D^1u_t|_{t=0}\|_{L^2}+\int_0^t[\|(\varrho, \theta)_\tau\|_{H^1}^2+
 \|(\varrho_\tau,u_\tau,\theta_\tau)\|_{H^1}^2 (\mathcal{E}^2+\|\varrho_\tau \|_{H^1}^2)]\mm{d}\tau.
\end{aligned}$$
Thus, using \eqref{detetim1} with $k=1$, \eqref{detetim2}  and Poincar\'e's inequality, we get \eqref{momentum2}.

Finally,  taking the time derivative in \eqref{ellipicequation},
 applying \eqref{ellioper} to the resulting identity, and making use of
\eqref{detetim1} with $k=1$ and \eqref{detetim2}, we obtain \eqref{utt}.
This completes the proof of Lemma \ref{lem:0304} \hfill$\Box$
\end{pf}

\subsection{Interior and boundary estimates }
To begin with, we establish the interior estimates of higher-order mass derivatives of $\bar{e}\varrho+\bar{\rho}\theta$.

\begin{lem}\label{lem:0305}    For $1\leq k\leq 3$, it holds that
$$\begin{aligned}
&\|\chi_0D^k(\varrho, u,\theta)(t)\|_{L^2}^2+\int_0^t\left(\|\chi_0  D^{k+1} u\|_{L^2}^2
+\left\|\chi_0D^k\frac{{d}}{{d}t}\left(\bar{e}\varrho+\bar{\rho}
\theta\right)\right\|_{L^2}^2\right)\mm{d}\tau \\
&\lesssim  \mathcal{E}_0^2+\int_0^t  \mathcal{R}\mm{d}\tau .
\end{aligned}$$
\end{lem}
\begin{pf}
Let $\chi_0$ be an arbitrary but fixed function in $C_0^\infty(\Omega)$. Then,
we can deduce from \eqref{masssimply}--\eqref{masssimply2} that
$$\begin{aligned}
&\int_0^t\int\left( \frac{a\chi_0^2\bar{e}}{\bar{\rho}}D^kL^\varrho D^k\varrho+
\chi_0^2D^k {L}^{u}\cdot D^k  u+\frac{\chi_0^2\bar{\rho}}{\bar{e}}D^kL^\theta D^k\theta\right)  \dxdt\\
&=\int_0^t\int\left(\frac{a\chi_0^2\bar{e}}{\bar{\rho}}D^kN^\varrho D^k\varrho+
\chi_0^2D^k {N}^{ u}\cdot D^k  u+\frac{\chi_0^2\bar{\rho}}{\bar{e}}D^kN^\theta D^k\theta \right)\dxdt,  \end{aligned}$$
which
yields that
\begin{equation}\label{esimttofinter}
\begin{aligned}
& \frac{1}{2}\int
\left(\frac{a\chi_0^2\bar{e}}{\bar{\rho}}|D^k \varrho|^2+\chi_0^2 {\bar{\rho}}| D^k u|^2+\frac{\chi_0^2\bar{\rho}}{\bar{e}}|D^k \theta|^2 \right)\mm{d}
\tau\\
&
+  \int_0^t\left(\mu\|\chi_0\nabla D^{k} u\|_{L^2}^2+\mu_0\|\chi_0  D^{k}\mm{div} u\|_{L^2}^2
 \right)\mm{d}\tau \\
&= \frac{1}{2}\int
\left(\frac{a\chi_0^2\bar{e}}{\bar{\rho}}|D^k \varrho(0)|^2+ \chi_0^2 \bar{\rho} |D^k u(0)|^2+\frac{\chi_0^2\bar{\rho}}{\bar{e}}|D^k \theta(0)|^2\right)\mm{d}t\\
&\quad -
2\int_0^t\int
\mu \chi_0(\nabla\chi_0\cdot\nabla D^k u ) \cdot D^k u\dxdt-
2\int_0^t\int
\mu \chi_0(\nabla\chi_0  D^k \mm{div}u ) D^k{ u}\dxdt
\\
&\quad -\int_0^t\int_\Omega
\chi_0^2[D^k(\bar{\rho} u_t)-\bar{\rho}D^k u_t]\cdot D^k u\dxdt-g\int_0^t\int \chi_0^2 D^k\varrho e_3\cdot D^ku\dxdt\\
&\quad-\int_0^t\int \left[\frac{a\chi_0^2\bar{e}}{\bar{\rho}}D^k(\bar{\rho}'u_3)D^k\varrho+\frac{\chi_0^2
\bar{\rho}}{\bar{e}}D^k(\bar{e}'u_3)D^k\theta\right]\dxdt\\
&\quad
+2a\int_0^t\int  \chi_0 D^k( \bar{e}\varrho+\bar{\rho}\theta)) \nabla\chi_0 \cdot D^ku
\dxdt\\
&\quad - \int_0^t\int
 \frac{a\chi_0^2\bar{e}}{\bar{\rho}} [D^k( \bar{\rho}\mm{div}u)-
 \bar{\rho}D^k\mm{div}u]   D^k\varrho
\dxdt\\ &
\quad -\int_0^t\int \frac{ a\chi_0^2\bar{\rho}}{\bar{e}} [D^k( \bar{e}\mm{div}u)-
 \bar{e}D^k\mm{div}u]    D^k\theta\dxdt\\
&\quad +\int_0^t\int\left(\frac{a\chi_0^2\bar{e}}{\bar{\rho}}D^kN^\varrho D^k\varrho+
\chi_0^2D^k {N}^{ u}\cdot D^k  u+\frac{\chi_0^2\bar{\rho}}{\bar{e}}D^kN^\theta D^k\theta \right)\dxdt:=L_1.
\end{aligned}
\end{equation}

Using \eqref{relationofdiv},  we can infer that
\begin{equation}\label{equaons1xdiv} \|\chi_0  D^{k}\mm{div} u\|_{L^2}^2
\geq \left \|\frac{\chi_0 }{\bar{\rho}\bar{e}+\bar{p}} D^{k}\frac{d}{dt}(\bar{e}\varrho+\bar{\rho}\theta)
\right\|_{L^2}^2
- \mathcal{R}
\end{equation}
On the other hand, exploiting the facts
$$-\int \frac{\chi_0^2\bar{\rho}}{\bar{e}}u\cdot \nabla D^k \theta D^k\theta\mm{d}x=
\frac{1}{2}\int | D^k \theta|^2\mm{div}\left(\frac{\chi_0^2\bar{\rho}}{\bar{e}} u\right)\mm{d}x \lesssim \mathcal{E}^3 $$
and $$\begin{aligned}
&-\int_0^t\int\chi_0^2 [\varrho  D^k u_t +aD^k \nabla (\varrho \theta)] \cdot D^ku\dxdt\\
&=-\frac{1}{2} \int\chi_0^2 \varrho  |D^k u|^2\mm{d}x+\frac{1}{2} \int\chi_0^2 \varrho_0  |D^k u_0|^2\mm{d}x
+ a\int_0^t\int D^k  (\varrho \theta)   \mm{div}(\chi_0^2 D^ku)\mm{d}x\\
& \quad\   +\frac{1}{2}\int_0^t\int\chi_0^2 (N^\varrho -\bar{\rho}' u_3-\bar{\rho}\mm{div}u)|D^k u|^2\dxdt\\
&\lesssim \mathcal{E}_0^2+ \int_0^t\mathcal{R}\mm{d}\tau+\left\|\frac{\varrho}{\bar{\rho}}\right\|_{L^\infty} \int\chi_0^2  {\bar{\rho}} |D^k u|^2(t)\mm{d}x,\end{aligned}$$
it's easy to derive that
\begin{equation}\label{equaons1xdiv2}
 \begin{aligned}L_1\leq &C\mathcal{E}_0^2+C\int_0^t  \mathcal{R}\mm{d}\tau
+C\left\|\frac{\varrho}{\bar{\rho}}\right\|_{L^\infty} \int\chi_0^2  {\bar{\rho}} |D^k u(t)|^2\mm{d}x \\
&+\frac{1}{2}\int_0^t\left(\mu\|\chi_0\nabla D^{k} u\|_{L^2}^2+\mu_0\|\chi_0  D^{k}\mm{div} u\|_{L^2}^2
 \right)\mm{d}\tau.
 \end{aligned}\end{equation}
Consequenetly, we get the desired conclusion  from \eqref{esimttofinter}, \eqref{equaons1xdiv} and \eqref{equaons1xdiv2}
immeditately.   \hfill$\Box$
\end{pf}

Next, let us establish the estimates near the boundary. Noting that $\partial\Omega$ is smooth,
 similarly to that in \cite{MANTIC481,MANTTP351},
we choose a finite number of bounded open sets $\{O^j\}_{j=1}^N$ in $\mathbb{R}^3$, such that $\cup_{j=1}^NO^j\supset \partial\Omega$.
In each open set $O^j$ we choose the local coordinates $y=(y_1,y_2,y_3)$ as follows:
 \begin{enumerate}
   \item[(1)] The surface $O^j\cap \partial\Omega$ is the image of a smooth vector function
   $z^j({y_1},y_2)=(z_1^j,z_2^j,z^j_3)({y_1},{y_2})$ (e.g., take the local geodesic polar coordinate), satisfying
    $|z_{y_1}^j|=1$, $z_{y_1}^j\cdot z_{y_2}^j=0$, and $|z_{y_2}^j|\geq \delta>0$, where $\delta$ is some positive constant independent of $1\leq j\leq N$.
      \item[(2)] Any $ x=(x_1,x_2,x_3)\in O^j$ is represented by
      \begin{equation}\label{transform}{x}_i=\omega_i(y):=
      y_3n_i(z^j({y_1},{y_2}))+z_i^j({y_1},{y_2})\mbox{ for }i=1,\ 2,\ 3,\end{equation} where $(n_1,n_2,n_3)(z^j({y_1},{y_2}))$ represents the internal unit normal vector
      at the point $z^j({y_1},{y_2})$ of the surface $\partial\Omega$.
 \end{enumerate}

For the simplicity of presentation, we omit the subscript $j$ in what follows. For $k=1,2$, we define the unit vectors
$$\tilde{e}_1=z_{y_1}\mbox{ and }\tilde{e}_2=z_{y_2}/|z_{y_2}|.$$
 An elementary calculation shows that the Jacobian $J$ of the transform \eqref{transform} is
\begin{equation}\label{jajobit}
J=\omega_{y_1}\times  \omega_{y_2}\cdot \tilde{n}=|z_{y_2}|+(\alpha|z_{y_2}|+\beta')y_3
 +(\alpha\beta'-\beta\alpha')y_3^2,\end{equation}
 where $\tilde{n}=(n_1,n_2,n_3)(z^j({y_1},{y_2})) $, $\alpha=\tilde{n}_{y_1}\cdot \tilde{e}_1$, $\beta=\tilde{n}_{y_1}\cdot \tilde{e}_2 $,  $\alpha'=\tilde{n}_{y_2}\cdot \tilde{e}_1$ and $\beta'=\tilde{n}_{y_2}\cdot \tilde{e}_2 $.
By \eqref{jajobit}, we find the transform \eqref{transform} is regular by choosing $y_3$ so small that $J\geq \delta/2$.
Therefore, the  inverse function of $\omega(y):=(\omega_1,\omega_2,\omega_3)(y)$ exits, and we denote it by $y=\omega^{-1}( x)$; moreover $({y_1},{y_2},y_3)_{x_i}( x)$  make sense and can be expressed by, using a straightforward calculation,
\begin{equation}\label{relations}\displaystyle
\left\{\begin{array}{l}
   \displaystyle \partial_{x_j}{y_1} =\frac{1}{J}( \omega_{y_2}\times  \omega_{y_3})_j=\frac{1}{J}(A\tilde{e}_j^1+B\tilde{e}_j^2)=:a_{1j}, \\[0.8em]
\displaystyle \partial_{x_j}{y_2} =\frac{1}{J}( \omega_{y_3}\times  \omega_{y_1})_j=\frac{1}{J}(C\tilde{e}_j^1+\tilde{D} \tilde{e}^2_j)=:a_{2j},  \\[0.8em]
 \displaystyle  \partial_{x_j}{y_3} =\frac{1}{J} (\omega_{y_1}\times  \omega_{y_2})_j=\tilde{n}_j=:a_{3j},
  \end{array}\right.
\end{equation}
where $A=|z_{y_2}|+\beta'y_3$, $B=-y_3\alpha'$, $C=-\beta y_3$, $\tilde{D}=1+\alpha y_3$,
\begin{equation}\label{jabjodbeiderelt}
J=A\tilde{D}-BC\geq \delta/2 \end{equation} and $\tilde{e}_j^m$
denotes the $j$-th component of $\tilde{e}_m$.
Obviously, \eqref{relations} gives
  $$\sum_{j=1}^3 a_{3j}^2=|\tilde{n}|^2=1,\   a_{1j}a_{3j}= a_{2j}a_{3j}=0,\ J^2=(AC+B\tilde{D})^2-(A^2+B^2)(C^2+\tilde{D}^2)$$ and
\begin{equation}\label{relationbetw}
 {\partial_{x_j}}=a_{kj} {\partial_{y_k}},
\end{equation}where we have
used the Einstein convention of summing over repeated indices.

Thus, in $O$, the three linear parts ${L}^{{\varrho}}$, ${L}^{{u}}_{\mm{new}}=(
{L}^{{u}^1}_{\mm{new}},{L}^{{u}^2}_{\mm{new}},{L}^{{u}^3}_{\mm{new}})$ and ${L}^{ {\theta}}$
in the local coordinates $({y_1},{y_2},y_3)$ read as follows.
\begin{equation*}\label{}\begin{aligned}
 &\displaystyle  {L}^{{\varrho}}=\tilde{L}^{\tilde{\varrho}}:=\tilde{\varrho}_t+
(a_{k 3}\partial_{y_k}\tilde{\rho})\tilde{u}_3+ {\tilde{\rho}}
  a_{kl} \tilde{u}^l_{y_k},\\
& \displaystyle {L}^{{u}^i}_{\mm{new}}=\tilde{L}^{\tilde{u}^i}_{\mm{new}}:=\tilde{\rho} \tilde{u}^i_t-\frac{\mu}{J^2}[(A^2+
 B^2) \tilde{u}^i_{{y_1}{y_1}}+2(AC+B\tilde{D}) \tilde{u}^i_{{y_1}{y_2}}+
 (C^2+\tilde{D}^2)\tilde{u}_{{y_2}{y_2}}^i+J^2\tilde{u}_{y_3y_3}^i ]
\\
& \displaystyle
\qquad +\mbox{ less two order terms of }\tilde{u}^i+g\tilde{\varrho} e_3+ a_{ki}\left(\frac{\mu_0}{\tilde{\rho}\tilde{e}+\tilde{p}}
\tilde{G}+a\tilde{e}\tilde{\varrho}+a\bar{\rho}
\tilde{\theta}\right)_{y_k} ,
 \\
  \displaystyle & {L}^{ {\theta}}=\tilde{L}^{ \tilde{\theta}}:= \tilde{\theta}_t+(a_{k 3}\partial_{y_k}\tilde{e})\tilde{u}_3+ {a\tilde{e}}  a_{kl}  \tilde{u}_{y_k}^l,
\end{aligned}\end{equation*}
where $(\tilde{\rho},\tilde{e},\tilde{p}):=(\bar{\rho},\bar{e},\bar{p})|_{x=\omega(y)}$,
$(\tilde{\varrho},\tilde{u},\tilde{\theta})(t,y)=(\varrho,u,\theta)|_{x=\omega(y)}$,
and
$$\tilde{G}:= ( \tilde{e} \tilde{\varrho}+\tilde{\rho} \tilde{\theta})_t+
u_l a_{kl} { ( \tilde{e} \tilde{\varrho}+\tilde{\rho} \tilde{\theta})}_{y_k}.$$
Similarly, we define that
 \begin{align*}   &\tilde{N}^{\tilde{\varrho}}:=
 {N}^{ {\varrho}}\;\;\mbox{ is written in the local coordinates} , \\
 &\tilde{N}^{\tilde{u}}_{\mm{new}}:= {N}^{ {u}}_{\mm{new}}\;\; \mbox{ is written in the local coordinates},\\
&\tilde{N}^{ \tilde{\theta}}:={N}^{ {\theta}}\;\;\mbox{ is written in the local coordinates}.
\end{align*}
With the notations above in hand, we can further rewrite the equations \eqref{masssimply}, \eqref{masssimply2} and \eqref{newform}
 in the local coordinates $y$ as follows:
\begin{equation}\label{newsystein}
   \tilde{L}^{\tilde{\varrho}} =\tilde{N}^{\tilde{\varrho}},
 \
    \tilde{L}^{{u}}_{\mm{new}}= \tilde{N}^{ \tilde{u}}_{\mm{new}},  \
\tilde{L}^{ \tilde{\theta}}= \tilde{N}^{ \tilde{\theta}},
 \end{equation}
with initial and boundary conditions
\begin{equation*}
(\tilde{\varrho}_0,\tilde{  u}_0,\tilde{\theta}_0) :=(\tilde{\varrho},\tilde{  u},\tilde{\theta})|_{t=0}=(\varrho_0,{ u}_0,\theta_0)|_{x= {\omega}(y)}\quad\mbox{in }{\tilde{\Omega}}
\end{equation*}
and
\begin{equation*}
\tilde{  u}(t,y)|_{\partial \tilde{\Omega}\cap \{y_3=0\}}={  0}\quad \mbox{ for any
}t>0,
\end{equation*}
where $\tilde{\Omega}=\{y~|~ y=w^{-1}(x),\ x\in O\cap \Omega\}$ and
$\tilde{L}^{\tilde{u} }_{\mm{new}}=( \tilde{L}^{\tilde{u}^1}_{\mm{new}},
\tilde{L}^{\tilde{u}^2}_{\mm{new}},\tilde{L}^{\tilde{u}^3}_{\mm{new}})$.

Let $\chi$ be an arbitrary but fixed function in $C_0^\infty(O )$ and $\tilde{\chi}:=\chi|_{x= {\omega}(y)}$. Obviously,
$\tilde{\chi} D^{k}_{{y_1} {y_2}} \tilde{u}=0$ on $\partial\tilde{\Omega}$.
Now, we control derivatives in the tangential directions.
\begin{lem}\label{lem:0306}   For $1\leq k\leq 3$, it holds that
$$\begin{aligned}
&\| \tilde{\chi} D^k_{{y_1}{y_2}}(  \tilde{\varrho},   \tilde{u}, \tilde{\theta})(t)\|_{L^2(\tilde{\Omega})}^2+\int_0^t\left(\| \tilde{\chi}D^{k}_{{y_1}{y_2}}D_{y}^1  \tilde{u}\|_{L^2( \tilde{\Omega})}^2
+\left\| \tilde{\chi} D^k_{{y_1}{y_2}}
 \tilde{G}\right\|_{L^2( \tilde{\Omega})}^2\right)\mm{d}\tau\\
& \lesssim  \mathcal{E}_0^2 +\int_0^t \mathcal{R}\mm{d}\tau.
\end{aligned}$$
\end{lem}
\begin{pf} To begin with, we deduce from
\eqref{newsystein} that
\begin{equation}\label{2012221433}\begin{aligned}
&\int_0^t\int_{\tilde{\Omega}}\left( \frac{a\tilde{\chi}^2\tilde{e}}{\tilde{\rho}}D^k_{{y_1}{y_2}}\tilde{L}^{\tilde{\varrho }} D^k_{{y_1}{y_2}}\tilde{\varrho}+
\tilde{\chi}^2D^k_{{y_1}{y_2}} \tilde{L}^{\tilde{ u}}_{\mm{new}}\cdot D^k_{{y_1}{y_2}}  \tilde{u}+\frac{\tilde{\chi}^2\tilde{\rho}}{\tilde{e}}D^k_{{y_1}{y_2}}\tilde{L}^{\tilde{\theta }} D^k_{{y_1}{y_2}}\tilde{\theta}
\right)\dydt\\
&=\int_0^t\int_{\tilde{\Omega}}\left(\frac{a\tilde{\chi}^2\tilde{e}}
{\tilde{\rho}}D^k_{{y_1}{y_2}}\tilde{N}^{\tilde{\varrho }} D^k_{{y_1}{y_2}}\tilde{\varrho}+
\tilde{\chi}^2D^k_{{y_1}{y_2}} \tilde{N}^{\tilde{ u}}_{\mm{new}}\cdot D^k_{{y_1}{y_2}}  \tilde{u}+\frac{\tilde{\chi}^2\tilde{\rho}}{\tilde{e}}
D^k_{{y_1}{y_2}}\tilde{N}^{\tilde{\theta}} D^k_{{y_1}{y_2}}\tilde{\theta } \right)\dydt.\end{aligned}
\end{equation}
In view of \eqref{relationofdiv}, we know that
$$\begin{aligned}
a_{kl}\tilde{u}^l_{y_k}=&
-\frac{1}{\tilde{\rho}\tilde{e}+\tilde{p}}\tilde{G}+
\frac{(a_{k3}\tilde{e}_{y_k}
 \tilde{\varrho}+a_{k3}\tilde{\rho}_{y_k}
 \tilde{\theta})\tilde{ u}_3+\tilde{\rho}(\tilde{N}^{\tilde{\theta}}+\tilde{u}_l a_{kl}\tilde{\theta}_{y_k})-a_{k3}(\tilde{\rho} \tilde{e})_{y_k}\tilde{u}_3- \tilde{e}\tilde{\varrho}
a_{kl}u_{y_k}^l}{\tilde{\rho}\tilde{e}+\tilde{p}}\\
=:&-\frac{1}{\tilde{\rho}\tilde{e}+\tilde{p}}\tilde{G}+L_2,
\end{aligned}$$
By the  formula of integration by parts, we have
\begin{equation}\label{20122214335}
\begin{aligned}&\int_{\tilde{\Omega}}\tilde{\chi}^2 D^k_{y_1 y_2}
\left[a_{ki}\left(\frac{\mu_0}{\tilde{\rho}\tilde{e}+\tilde{p}}\tilde{G}\right)_{y_k}
\right]D^k_{y_1 y_2}
\tilde{u}_i\mm{d}y\\
=&
\int_{\tilde{\Omega}}
 \frac{\mu_0\tilde{\chi}^2}{(\tilde{\rho}\tilde{e}+\tilde{p})^2}|D^k_{y_1 y_2}\tilde{G}|^2\mm{d}y+
\int_{\tilde{\Omega}}(\mbox{less than } k+1 \mbox{ order terms of }\tilde{G})
 D^{k}_{y_1y_2}L_2\mm{d}y\\
&+
\int_{\tilde{\Omega}}(\mbox{less than } k+1 \mbox{ order terms of }\tilde{G}) (\mbox{less than } k \mbox{ order terms of }\tilde{G})
  \mm{d}y
\\
&+
\int_{\tilde{\Omega}}(\mbox{less than } k  \mbox{ order terms of }\tilde{G})
(\mbox{less than } k+2  \mbox{ order terms of }\tilde{u} ) \mm{d}y\\
=:&\int_{\tilde{\Omega}}
 \frac{\mu_0\tilde{\chi}^2}{(\tilde{\rho}\tilde{e}+\tilde{p})^2}|D^k_{y_1 y_2}\tilde{G}|^2\mm{d}y
+L_3.
\end{aligned}
\end{equation}
Thus, similarly to \eqref{esimttofinter}, we can infer from the  two equalities \eqref{2012221433} and \eqref{20122214335}  that
\begin{equation}\label{newesimtefro}
\begin{aligned}
& \frac{1}{2}\int_{\tilde{\Omega}}
\left(\frac{a\tilde{\chi}^2\tilde{e}}{\tilde{\rho}}|D^k_{{y_1}{y_2}} \tilde{\varrho}|^2+ \tilde{\chi}^2\tilde{\rho}|D^k_{{y_1}{y_2}} \tilde{u}|^2+\frac{\tilde{\chi}^2\tilde{\rho}}{\tilde{e}}|D^k_{{y_1}{y_2}} \tilde{\theta}|^2 \right)(t)\mm{d}
y\\
&+\int_0^t\int_{\tilde{\Omega}}
 \frac{\mu_0\tilde{\chi}^2}{(\tilde{\rho}\tilde{e}+\tilde{p})^2}|D^k_{y_1 y_2}\tilde{G}|^2\dydt+  \int_0^t\int_{\tilde{\Omega}}
\frac{\mu\tilde{\chi}^2}{J^2}\sum_{i=1}^3[(A^2+
 B^2) | D^k_{{y_1}{y_2}}\tilde{u}^i_{{y_1} }|^2\\
&\quad   +2(AC+B\tilde{D}) D^k_{{y_1}{y_2}}\tilde{u}^i_{{y_1}}D^k_{{y_1}{y_2}}\tilde{u}^i_{y_2}+
(C^2+\tilde{D}^2)| D^k_{{y_1}{y_2}}\tilde{u}_{{y_2} }^i|^2+J^2| D^k_{{y_1}{y_2}}\tilde{u}_{y_3}^i|^2 ]
\mm{d}y\mm{d}\tau \\
&= \frac{1}{2}\int_{\tilde{\Omega}}
\left(\frac{a\tilde{\chi}^2\tilde{e}}{\tilde{\rho}}|D^k_{{y_1}{y_2}} \tilde{\varrho}_0|^2+  \tilde{\chi}^2\tilde{\rho}|D^k_{{y_1}{y_2}} \tilde{u}_0|^2+\frac{\tilde{\chi}^2\tilde{\rho}}{\tilde{e}}|D^k_{{y_1}{y_2}} \tilde{\theta}_0|^2\right)\mm{d}t\\
&\quad -\int_0^t\int_{\tilde{\Omega}}
\tilde{\chi}^2[D^k_{{y_1}{y_2}}(\tilde{\rho}\tilde{ u}_t)-\tilde{\rho} D^k_{{y_1}{y_2}} \tilde{u}_t]\cdot D^k_{{y_1}{y_2}} \tilde{u}\dydt\\
&\quad +\sum_{i=1}^3\int_0^t\int
(\mbox{less than } k+1 \mbox{  order terms of } \tilde{\varrho}\mbox{ and } \tilde{\theta}) \\
&\qquad\qquad\qquad \ ( \mbox{less than } k+1 \mbox{ order of }\tilde{u}_i)\dydt \\
&\quad +\sum_{j=1}^3\sum_{i=1}^3\int_0^t\int
(\mbox{less than } k+2 \mbox{  order terms of } \tilde{u}_i) D^k_{{y_1}{y_2}} \tilde{u}_j\dydt
\\
&\quad +\sum_{j=1}^2\sum_{i=1}^3\int_0^t\int
(\mbox{less than } k+1 \mbox{  order terms of } \tilde{u}_i) D^k_{{y_1}{y_2}} \tilde{h}_j \dydt\\
&\quad
+\int_0^t\int_{\tilde{\Omega}}\left(\frac{a\tilde{\chi}^2\tilde{e}}{\tilde{\rho}}
D^k_{{y_1}{y_2}}\tilde{N}^{\tilde{\varrho}} D^k_{{y_1}{y_2}}\tilde{\varrho}+
\tilde{\chi}^2D^k_{{y_1}{y_2}} \tilde{N}^{\tilde{ u}}_{\mm{new}}\cdot D^k_{{y_1}{y_2}}  \tilde{u}+\frac{\tilde{\chi}^2\tilde{\rho}}{\tilde{e}}D^k_{{y_1}{y_2}}
\tilde{N}_{\mm{new}}^{\tilde{\theta}} D^k_{{y_1}{y_2}}\tilde{\theta} \right)\dydt\\
&\quad +L_3:=L_4,
\end{aligned}
\end{equation}where we have defined $\tilde{h}_1=\tilde{\varrho}$  and $\tilde{h}_2=\tilde{\theta}$ for simplicity.

By virtue of \eqref{jabjodbeiderelt}, the matrix
$$\left(
  \begin{array}{cc}
    A^2+B^2 &  AC+B\tilde{D} \\
    AC+B\tilde{D}  & C^2+\tilde{D}^2, \\
  \end{array}
\right)$$
is strictly positive-defined, there thus exists a positive constant $\underline{\mu}$, such that
\begin{equation}\label{201512427}\begin{aligned}&\underline{\mu} D^{k}_{{y_1}{y_2}}| \tilde{u}_{y_j}^i|^2
\leq  (A^2+
 B^2) |  D^k_{{y_1}{y_2}}\tilde{u}^i_{{y_1} }|^2+2(AC+B\tilde{D}) D^k_{{y_1}{y_2}}\tilde{u}^i_{{y_1}}D^k_{{y_1}{y_2}}\tilde{u}^i_{y_2}\\
&\qquad \qquad \qquad +
 (C^2+\tilde{D}^2)| D^k_{{y_1}{y_2}}\tilde{u}_{{y_2} }^i|^2+J^2| D^k_{{y_1}{y_2}}\tilde{u}_{y_3}^i|^2.
\end{aligned}\end{equation}

On the other hand, similarly to \eqref{equaons1xdiv2}, it is easy to see that
\begin{equation}
\label{2015121425}\begin{aligned}
L_4\leq &C \tilde{\mathcal{E}}_0^2 +C\int_0^t \tilde{\mathcal{R}}\mm{d}\tau+\frac{\underline{\mu}}{2}
\int_0^t \int_{\tilde{\Omega}}\tilde{\chi}^2|D^{k}_{{y_1}{y_2}}D_{y}^1 \tilde{u} |^2 \dydt\\
&\  +\frac{1}{2}\int_0^t\int_{\tilde{\Omega}}
 \frac{\mu_0\tilde{\chi}^2}{(\tilde{\rho}\tilde{e}+\tilde{p})^2}|D^k_{y_1 y_2}\tilde{G}|^2\dydt
 + C\left\|\frac{\tilde{\varrho}}{\tilde{\rho}}\right\|_{L^\infty(\tilde{\Omega})} \int_{\tilde{\Omega}}
  \tilde{\chi}^2\tilde{\rho}|D^k_{{y_1}{y_2}} \tilde{u}|^2(t)  \mm{d}
y,
\end{aligned}
\end{equation}
where we have defined
$\tilde{\mathcal{E}}(t):=\tilde{\mathcal{E}}(\tilde{\varrho},\tilde{u},
\tilde{\theta})(t):=\|(\tilde{\varrho}, \tilde{u},\tilde{\theta})(t)\|_{H^3
(\tilde{\Omega})}$,  $\tilde{\mathcal{E}}_0=\tilde{\mathcal{E}}(0)$ and
$$ \tilde{\mathcal{R}}:=\left\|\left(\tilde{\varrho},
\tilde{\theta},\tilde{u}_t,\frac{{d}}{{d}t}
\left(\bar{e}\tilde{\varrho}+\bar{\rho}
\tilde{\theta}\right)\right)
\right\|_{H^2(\tilde{\Omega})}^2+
\tilde{\mathcal{E}} (\|\tilde{u} \|_{H^3(\tilde{\Omega})}+
\|\tilde{u} \|_{H^4(\tilde{\Omega})}^2
+\tilde{\mathcal{E}}^2 ).$$
 Utilizing \eqref{2015121425} and \eqref{201512427}, we deduce from \eqref{newesimtefro} that
\begin{equation}\label{2015124289}\begin{aligned}
&\| \tilde{\chi} D^k_{{y_1}{y_2}}(  \tilde{\varrho},   \tilde{u}, \tilde{\theta})(t)\|_{L^2(\tilde{\Omega})}^2+\int_0^t\left(\| \tilde{\chi}D^{k}_{{y_1}{y_2}}D_{y}^1 \tilde{u}\|_{L^2( \tilde{\Omega})}^2
+\left\| \tilde{\chi} D^k_{{y_1}{y_2}}\tilde{G}
 \right\|_{L^2( \tilde{\Omega})}^2\right)\mm{d}\tau\\
& \lesssim  \tilde{\mathcal{E}}_0^2 +\int_0^t \tilde{\mathcal{R}}\mm{d}\tau,
\end{aligned}\end{equation}

Finally, using \eqref{transform} and \eqref{jajobit}, we
can obtain the following estimate by transformation of coordinates
\begin{equation}
\label{trnasesimate}
\|f\|_{H^k(\tilde{\Omega})}\lesssim \|f\|_{H^k }.
\end{equation}
In particular, we have $\tilde{\mathcal{E}}_0\lesssim  {\mathcal{E}}_0$  and
$\tilde{\mathcal{R}}\lesssim  \mathcal{R} $. Consequently, we obtain the desired conclusion from \eqref{2015124289}.
\hfill$\Box$
\end{pf}

Next, we turn to the estimate of derivatives in the normal directions.
\begin{lem}\label{lem:0307}
 For $0\leq  k+l\leq  2$, it holds that
$$\begin{aligned}
&\| \tilde{\chi} D_{{y_1}{y_2}}^k D_{y_3}^{l+1}(\tilde{e} \tilde{\rho}+\tilde{\rho} \tilde{\theta})(t)\|_{L^2( \tilde{\Omega})}^2\\
&\quad +\int_0^t\left\{\| \tilde{\chi} D_{{y_1}{y_2}}^k D_{y_3}^{l+1}(\tilde{e} \tilde{\varrho}+\tilde{\rho} \tilde{\theta})\|_{L^2( \tilde{\Omega})}^2+\left\|\chi D_{{y_1}{y_2}}^k D_{y_3}^{l+1}\tilde{G}\right\|_{L^2(\tilde{\Omega})}^2\right\}
\mm{d}\tau
\\
&\lesssim \|(\varrho_0,\theta_0)\|_{H^{3}}^2+\int_0^t(\|D_{{y_1}{y_2}}^{k+1}D_{y_3}^lD_{
y}^1 \tilde{u}\|_{L^2( \tilde{\Omega})}^2+ \mathcal{R})\mm{d}\tau.
\end{aligned}$$
\end{lem}
\begin{pf} Recalling that $\sum_{i=1}^3\tilde{n}_ia_{ki}=\sum_{i=1}^3a_{3i}a_{ki}=0$, $\sum_{i=1}^3a_{3i}^2=1$ and the relation \eqref{relations}, we use the equations
$D_{y_3}^1(\tilde{e}\tilde{L}^{ \tilde{\varrho}}+\tilde{\rho}\tilde{L}^{ \tilde{\theta}}-
\tilde{e} \tilde{N}^{ \tilde{\varrho}}-\tilde{\rho} \tilde{N}^{ \tilde{\theta}})=0$ and $\tilde{n}\cdot(\tilde{L}^{ \tilde{u}}_{\mm{new}}
-\tilde{N}^{ \tilde{u}}_{\mm{new}})=0$, to arrive at
$$\begin{aligned}
&\tilde{G}_{y_3}+\frac{\tilde{\rho}\tilde{e}+\tilde{p}}{J}[(A\tilde{e}_1+B\tilde{e}_2)\cdot  \tilde{u}_{{y_3}{y_1}}+
  (C\tilde{e}_1+\tilde{D}\tilde{e}_2)\cdot \tilde{ u}_{{y_3}{y_2}}+J \tilde{n}\cdot  \tilde{u}_{{y_3}{y_3}}]\\
  &\quad +\mbox{less than two order terms of } \tilde{u}\\
 &=[\tilde{\rho}( \tilde{N}^{ \tilde{\theta}}+ \tilde{u}_{i}a_{ki}
 \tilde{\theta}_{y_k})-\tilde{e} \tilde{\varrho}a_{ki} \tilde{u}^i_{y_k}
  +(a_{k3}\tilde{e}_{y_k} \tilde{\varrho}-a_{k3}\tilde{\rho}_{y_k} \tilde{\theta} ) \tilde{u}_3]_{y_3}
  \end{aligned}$$
and
\begin{equation}\label{moetnromal}\begin{aligned}
&\displaystyle  \tilde{\rho} \tilde{n}\cdot \tilde{u}_t-
 \frac{\mu  \tilde{n}}{J^2}\cdot[(A^2+B^2) \tilde{u}_{{y_1}{y_1}}+2(AC+B\tilde{D}) \tilde{u}_{{y_1}{y_2}}+
 (C^2+\tilde{D}^2) \tilde{u}_{{y_2}{y_2}}+J^2 \tilde{u}_{{y_3}{y_3}} ] \\
& \displaystyle +\mbox{less than two order terms of } \tilde{u}
+g \tilde{\varrho} e_3\cdot  \tilde{n} +\left[\frac{\mu_0}{\tilde{\rho}\tilde{e}+\tilde{p}}\tilde{G}+
a\tilde{e} \tilde{\varrho}+ a\tilde{\rho} \tilde{\theta}\right]_{y_3} \\
& =  \tilde{n}\cdot\tilde{N}^{ \tilde{u}}_{\mm{new}}   \end{aligned}
\end{equation}
Eliminating $\mu  \tilde{n}\cdot  \tilde{u}_{rr}$ from \eqref{moetnromal}, we get
\begin{equation}\label{moetnromal2}
\begin{aligned}
\displaystyle &\left[\frac{(\mu+\mu_0)}{\tilde{\rho}\tilde{e}+\tilde{p}}
\tilde{G}+
a\tilde{e} \tilde{\varrho}+a\tilde{\rho} \tilde{\theta}\right]_{y_3}
=- \tilde{\rho} \tilde{n}\cdot \tilde{u}_t+
 \frac{\mu \tilde{n}}{J^2}\cdot [(A^2+B^2) \tilde{u}_{{y_1}{y_1}}+2(AC+B\tilde{D}) \tilde{u}_{{y_1}{y_2}}   \\
&\quad    +  (C^2+\tilde{D}^2) \tilde{u}_{{y_2}{y_2}} ]-\frac{\mu}{J}[(A\tilde{e}_1+B\tilde{e}_2)\cdot  \tilde{u}_{{y_3}{y_1}}+
  (C\tilde{e}_1+ \tilde{D}\tilde{e}_2)\cdot  \tilde{u}_{{y_3}{y_2}}]     \\
& \quad    +\mbox{less than two } \mbox{order terms of } \tilde{u} \displaystyle-g \tilde{\varrho} e_3\cdot
\tilde{n}+ \tilde{n}\cdot\tilde{N}^{ \tilde{u}}_{\mm{new}}  \\
&\quad +\frac{\mu}{\tilde{\rho}\tilde{e}+\tilde{p}}[\tilde{\rho}( \tilde{N}^{ \tilde{\theta}}+ \tilde{u}_{i}a_{ki}
 \tilde{\theta}_{y_k})
  +(a_{k3}\tilde{e}_{y_k} \tilde{\varrho}+a_{k3}\tilde{\rho}_{y_k} \tilde{\theta} ) \tilde{u}_3-\tilde{e} \tilde{\varrho}a_{ki} \tilde{u}^i_{y_k}]_{y_3}:=L_5.
\end{aligned}
\end{equation}

If we apply $D_{{y_1}{y_2}}^kD_{y_3}^l$ ($k+l=0$, $1$, $2$) to
 \eqref{moetnromal2}, multiply then by
$\tilde{\chi}^2[ D_{{y_1}{y_2}}^kD_{y_3}^l\tilde{G}_{y_3}+ D_{{y_1}{y_2}}^{k}D_{y_3}^{l+1}(\tilde{e} \tilde{\varrho} +\tilde{\rho} \tilde{\theta})]$
and integrate them, we get
\begin{equation*}
\begin{aligned}
\displaystyle &\frac{1}{2}\int_{\tilde{\Omega}}
 \left[\frac{(\mu+\mu_0)}{\tilde{\rho}
\tilde{e}+\tilde{p}}+a\right]\tilde{\chi}^2| D_{{y_1}{y_2}}^k D_{y_3}^{l+1}(\tilde{e} \tilde{\varrho}+\tilde{\rho} \tilde{\theta})(t)|^2\mm{d}y\\
&+\int_0^t\int_{\tilde{\Omega}}\left(a|\tilde{\chi} D_{{y_1}{y_2}}^k D_{y_3}^{l+1}( \tilde{e} \tilde{\varrho}+ \tilde{\rho} \tilde{\theta})|^2 + \frac{ \mu+\mu_0 }{\tilde{\rho}
\tilde{e}+\tilde{p}}| \tilde{\chi} D_{{y_1}{y_2}}^k D_{y_3}^{l+1}\tilde{G}|^2\right)
\dydt\\
&=\frac{1}{2}\int_{\tilde{\Omega}}
 \left[\frac{(\mu+\mu_0)}{\tilde{\rho}
\tilde{e}+\tilde{p}}+a\right]\tilde{\chi}^2| D_{{y_1}{y_2}}^k D_{y_3}^{l+1}(\bar{e} \tilde{\varrho}_0+\bar{\rho} \tilde{\theta}_0) |^2\mm{d}y\\
&\quad +\int_0^t\int_{\tilde{\Omega}} (\mbox{less than }k+l+1\mbox{ order terms of }\tilde{G} + D_{{y_1}{y_2}}^kD_{y_3}^l L_5 )\\
&\qquad \qquad \quad \,\tilde{\chi}^2[ D_{{y_1}{y_2}}^kD_{y_3}^l\tilde{G}_{y_3}+ D_{{y_1}{y_2}}^kD_{y_3}^{l+1}(\tilde{e} \tilde{\varrho} +\tilde{\rho} \tilde{\theta}) ]\dydt\\
&\quad -\int_0^t\int_{\tilde{\Omega}} \tilde{\chi}^2 \{D_{{y_1}{y_2}}^kD_{y_3}^{l+1}
 [u_m a_{jm}  ( \tilde{e} \tilde{\varrho}+\tilde{\rho} \tilde{\theta})_{y_j}]-
u_m a_{jm} D_{{y_1}{y_2}}^kD_{y_3}^{l+1} ( \tilde{e} \tilde{\varrho}+\tilde{\rho} \tilde{\theta})_{y_j}\} \\
&\qquad \qquad \quad \, \left(\frac{ \mu+\mu_0 }{\tilde{\rho}\tilde{e}+\tilde{p}} +a\right)D_{{y_1}{y_2}}^kD_{y_3}^{l+1}(   \tilde{e} \tilde{\varrho}+a\tilde{\rho} \tilde{\theta}) \dydt\\
&\quad \,+\frac{1}{2}\int_0^t\int_{\tilde{\Omega}}\left[ \tilde{\chi}^2 u_m a_{jm} \left(a+
\frac{ \mu+\mu_0 }{ \tilde{\rho}\tilde{e}+\tilde{p}}\right)\right]_{y_i} |D_{{y_1}{y_2}}^kD_{y_3}^{l+1}  ( \tilde{e} \tilde{\varrho}+\tilde{\rho} \tilde{\theta}) ]|^2  \dydt\\
&:=L_6.
\end{aligned}
\end{equation*}
Noting that $ H^2\hookrightarrow L^\infty$ and $\mathcal{E}\leq \delta_1^0$ is sufficiently small, thus it's easy to estimate that
$$\begin{aligned}
L_6\lesssim &\|(\tilde{\varrho}_0,\tilde{\theta}_0)\|_{H^{3}(\tilde{\Omega})}^2+
\int_0^t(\|D_{{y_1}{y_2}}^{k+1}D_{y_3}^lD \tilde{u}\|_{L^2( \tilde{\Omega})}^2+ \tilde{\mathcal{R}})\mm{d}\tau\\
&+\frac{1}{2}\int_0^t\int_{\tilde{\Omega}}\left(a|\tilde{\chi} D_{{y_1}{y_2}}^k D_{y_3}^{l+1}(\tilde{e} \tilde{\varrho}+\tilde{\rho} \tilde{\theta})|^2 + \frac{ \mu+\mu_0 }{\tilde{\rho}
\tilde{e}+\tilde{p}}|\tilde{\chi} D_{{y_1}{y_2}}^k D_{y_3}^{l+1}\tilde{G}|^2\right) \dydt,
\end{aligned}$$
and
$$\begin{aligned}
&\| \tilde{\chi} D_{{y_1}{y_2}}^k D_{y_3}^{l+1}(\tilde{e} \tilde{\rho}+\tilde{\rho} \tilde{\theta})(t)\|_{L^2( \tilde{\Omega})}^2\\
&\quad +\int_0^t\left\{\| \tilde{\chi} D_{{y_1}{y_2}}^k D_{y_3}^{l+1}(\tilde{e} \tilde{\varrho}+\tilde{\rho} \tilde{\theta})\|_{L^2( \tilde{\Omega})}^2+\left\|\tilde{\chi} D_{{y_1}{y_2}}^k D_{y_3}^{l+1}\tilde{G}\right\|_{L^2(\tilde{\Omega})}^2\right\}
\mm{d}\tau
\\
&\lesssim \|(\tilde{\varrho}_0,\tilde{\theta}_0)\|_{H^{3}(\tilde{\Omega}}^2+
\int_0^t(\|D_{{y_1}{y_2}}^{k+1}D_{y_3}^lD \tilde{u}\|_{L^2( \tilde{\Omega})}^2+ \tilde{\mathcal{R}})\mm{d}\tau,
\end{aligned}$$
which, together with \eqref{trnasesimate}, yields the desired conclusion.
\hfill$\Box$
\end{pf}

Finally, we introduce the following lemma \cite[Lemma 5.8]{MANTTP351} on the stationary Stokes equations to
get the estimates on the tangential derivatives of both $u$ and $\bar{e}\varrho+\bar{\rho}\theta$.
\begin{lem}\label{lem:0302}
Consider the Stokes problem
$$\left\{\begin{array}{l}
   -\mu\Delta u+a\nabla \sigma=g,  \\
   \mm{div}u=f, \\
  u|_{\partial \Omega}=0,
\end{array}\right.$$
where $f\in H^{k+1}$ and $g\in H^k$ ($k\geq 0$). Then the above problem
  has a solution $(\sigma,u)\in H^{k+1}\times H^{k+2}\cap H_0^{1}$ which is unique modulo a constant of integration
  for $\sigma$. Moreover, this solution satisfies
$$\|u\|_{H^{k+2}}^2+\|D\sigma\|_{H^k}^2\lesssim \|f\|_{H^{k+1}}^2+\|g\|_{H^k}^2.$$
\end{lem}
Now we rewrite the perturbed equations as the  Stokes problem:
\begin{equation*}\label{stokesu}\left\{\begin{array}{ll}
-\mu \Delta  u+a\nabla(\bar{e} {\varrho}  +\bar{\rho}\theta ) ={N}^u-\bar{\rho} u_t-g\varrho  e_3+\mu_0\nabla \mm{div} u,
\\
\displaystyle ({\bar{\rho}\bar{e}+\bar{p}})\mm{div} u
 =\bar{\rho}(N^\theta+u\cdot\nabla \theta)-{\bar{e}\varrho\mm{div}u+(\bar{e}'\varrho+\bar{\rho}'\theta)u_3-
 \frac{d}{dt}(\bar{e}\varrho+\bar{\rho}\theta)}-(\bar{\rho}\bar{e})'{u}_3,
 \\
u|_{\partial\Omega}=0.
 \end{array}\right.\end{equation*}
which, together with \eqref{transform}, implies that the following  Stokes problem of $(\tilde{\chi} D^k_{{y_1}{y_2}}\tilde{u}, \tilde{\chi} D^k_{{y_1}{y_2}}(\bar{e}
 \tilde{\varrho}  +\bar{\rho}\tilde{\theta} )) |_{y=\omega^{-1}(x)}$:
  \begin{equation*}\label{}\left\{\begin{array}{ll}
-\mu \Delta  [(\tilde{\chi} D^k_{{y_1}{y_2}}\tilde{u})|_{y=\omega^{-1}(x)}]+a\nabla\{[(\tilde{\chi} D^k_{{y_1}{y_2}}(\tilde{e}
 \tilde{\varrho}  +\tilde{\rho}\tilde{\theta} )]|_{y=\omega^{-1}(x)}\} \\
\; =\mbox{less than } k+2 \mbox{ order terms of }u \\
  \quad +\mbox{less than } k+1 \mbox{ order terms of }\varrho,\ \theta\mbox{ and } ({N}^u-
\bar{\rho} u_t-g\varrho  e_3 )\\
\quad  + \mu_0(\tilde{\chi}D^k_{{y_1}{y_2}} \tilde{U})|_{y=\omega^{-1}(x)}\ +\mbox{ less than } k  \mbox{ order terms of }   \nabla \mm{div} u :=L_7,\\[2mm]
  \displaystyle  \mm{div} [(\tilde{\chi}  D^k_{{y_1}{y_2}}\tilde{u})|_{y=\omega^{-1}(x)}]
 = \\
 \quad  -  \{\tilde{\chi} D^k_{{y_1}{y_2}} [ \frac{d}{dt}(\tilde{e}\tilde{\varrho}+\tilde{\rho}\tilde{\theta})
({\tilde{\rho}\tilde{e}+\tilde{p}})^{-1}]\}_{y=\omega^{-1}(x)}
+\mbox{less than } k  \mbox{ order terms of }\\
\quad \frac{d}{dt}(\bar{e}\varrho+\bar{\rho}\theta)
 + \mbox{less than } k+1 \mbox{ order  terms of } u\mbox{ and }\\
\quad \{ [\bar{\rho}(N^\theta+u\cdot\nabla \theta)
-\bar{e}\varrho\,\mm{div}u+(\bar{e}'\varrho+\bar{\rho}'\theta)u_3
-(\bar{\rho}\bar{e})'{u}_3 ]
 ({\bar{\rho}\bar{e}+\bar{p}})^{-1} \}:=L_8, \\[1mm]
[ ( \tilde{\chi} D_{{y_1}{y_2}}^k\tilde{u})|_{y=\omega^{-1}(x)}]|_{\partial\Omega}=0,
 \end{array}\right.\end{equation*}where $\tilde{U}:= \nabla \mm{div} u  $ is written in local  coordinates.
 Applying Lemma \ref{lem:0302} to the above problem, we obtain
$$\begin{aligned}
&\| D^{2+l}_{x}[(\tilde{\chi} D_{{y_1}{y_2}}^k  \tilde{u})|_{y=\omega^{-1}(x)}]\|_{L^2 }^2+
\| D^{1+l}_{x}[\tilde{\chi} D_{{y_1}{y_2}}^k (\tilde{e} \tilde{\varrho}+\tilde{\rho} \tilde{\theta})]|_{y=\omega^{-1}(x)}\|_{L^2 }^2\\
&\lesssim \| L_8\|_{H^{l+1}}^2+\| L_7 \|_{H^{l}}^2\mbox{ for }0\leq l+k\leq 2,
\end{aligned}$$
which, combined with the estimates
 $$
 \| L_8\|_{H^{l+1}}^2+\| L_7 \|_{H^{l}}^2\lesssim \left\| {\chi} (x)D^{1+l}_x ( D_{{y_1}{y_2}}^{k}\tilde{G}) |_{y=\omega^{-1}(x)}\right\|_{L^2 }^2
+ \left\| {\chi}(x)D^{l}_x(D^k_{{y_1}{y_2}} \tilde{U})|_{y=\omega^{-1}(x)}\right\|_{L^2 }^2+\mathcal{R},
$$
and  $$\left\|{\chi}(x)D^{l}_x( D^k_{{y_1}{y_2}} \tilde{U})|_{y=\omega^{-1}(x)}\right\|_{L^2 }^2\lesssim  \left\| {\chi} (x)D^{1+l}_x ( D_{{y_1}{y_2}}^{k}\tilde{G}) |_{y=\omega^{-1}(x)}\right\|_{L^2 }^2+\mathcal{R}, $$
yields
$$\begin{aligned}
&\|  {\chi}(x)  D^{2+l}_{x}(D_{{y_1}{y_2}}^k  \tilde{u})|_{y=\omega^{-1}(x)}]\|_{L^2 }^2+
\|  {\chi}(x)  D^{1+l}_{x} [D_{{y_1}{y_2}}^k (\tilde{e} \tilde{\varrho}+\tilde{\rho} \tilde{\theta})]|_{y=\omega^{-1}(x)}\|_{L^2 }^2
\\
&\lesssim \left\| {\chi} (x)D^{1+l}_x ( D_{{y_1}{y_2}}^{k}\tilde{G}) |_{y=\omega^{-1}(x)}\right\|_{L^2 }^2
 +\mathcal{R} .
\end{aligned}$$
Consequently, by the transformation of coordinates  \eqref{transform},
 we easily obtain the following estimate:
\begin{lem}\label{lem:0308}  For $0\leq l+k\leq 2$, we have
 $$\|\tilde{\chi}  D^{2+l}_yD_{{y_1}{y_2}}^k  \tilde{u}\|_{L^2( \tilde{\Omega})}^2+
\|\tilde{\chi}  D^{1+l}_yD_{{y_1}{y_2}}^k (\tilde{e} \tilde{\varrho}+\tilde{\rho} \tilde{\theta})\|_{L^2( \tilde{\Omega})}^2
\lesssim \left\|\tilde{\chi} D^{1+l}_y D_{{y_1}{y_2}}^{k}\tilde{G}\right\|_{L^2( \tilde{\Omega})}^2+\mathcal{R}. $$
\end{lem}

Now, we are able to establish the desired boundary estimate.
\begin{lem}\label{lem:0309} For $ 0\leq  k\leq  3$, it holds that$$ \sum_{k=0}^3\int_0^t\left\{\| {\chi}
D^{k+1}   {u}\|_{L^2 }^2+\left\|\chi D^k\left[\frac{d}{dt}(\bar{e} {\varrho}+\bar{\rho} {\theta})\right]\right\|_{L^2 }^2\right\}\mm{d}\tau
\lesssim \mathcal{E}_0^2 +\int_0^t\mathcal{R}\mm{d}\tau.$$
\end{lem}
\begin{pf}
It  suffices to show
\begin{equation}\label{sumefork=1} \int_0^t\left\{\| \tilde{\chi}
D^{k}_yD^1_y   \tilde{u}\|_{L^2(\tilde{\Omega}) }^2+\left\|\tilde{\chi} D^k_y\tilde{G}\right\|_{L^2
(\tilde{\Omega}) }^2\right\}\mm{d}\tau
\lesssim \mathcal{E}_0^2 +\int_0^t\mathcal{R}\mm{d}\tau,
\end{equation}
from which we can immediately get the desired conclusion by transformation of coordinates.

We consider $D_y^{k}=D_{y_3}^{l}D_{y_1y_2}^{k-l} $ for $1\leq l\leq 3$. Then,
$$\begin{aligned} & \int_0^t\left\{\|\tilde{\chi}
D_{y_3}^{l}D_{y_1y_2}^{k-l} D^{1 }_{y} \tilde{u}\|_{L^2( \tilde{\Omega})}^2
+\left\|\tilde{\chi} D_{y_3}^{l}D_{y_1y_2}^{k-l}\tilde{G}\right\|_{L^2( \tilde{\Omega})}^2\right\}\mm{d}\tau\\
&\lesssim  \int_0^t\left[\left\|\tilde{\chi} D^{l}_y D_{{y_1}{y_2}}^{k-l}
\tilde{G}\right\|_{L^2( \tilde{\Omega})}^2+\mathcal{R}\right]\mm{d}\tau\\
&\lesssim \|(\varrho_0,\theta_0)\|_{H^{3}}^2+ \int_0^t  \mathcal{R}\mm{d}\tau
+\sum_{0\leq m\leq l-1 } \int_0^t\|\tilde{\chi}D_{y_3}^{m}D_{{y_1}{y_2}}^{k-m}D_{
y}^1 \tilde{u}\|_{L^2( \tilde{\Omega})}^2 \mm{d}\tau,
\end{aligned}$$
where we have used  Lemmas \ref{lem:0308} and \ref{lem:0307} for the first and  second inequalities respectively.
Obviously, if $m\geq 1$, we can continue to repeat the  process above. Thus we conclude that
$$\begin{aligned}  & \int_0^t\left\{\|\tilde{\chi}
D^{k+1}_{y} \tilde{u}\|_{L^2( \tilde{\Omega})}^2
+\left\|\tilde{\chi}  D_y^{k}\tilde{G}\right\|_{L^2( \tilde{\Omega})}^2\right\}\mm{d}\tau \\
&\lesssim \|(\varrho_0,\theta_0)\|_{H^{3}}^2+ \int_0^t  \mathcal{R}\mm{d}\tau+\int_0^t\|\tilde{\chi} D_{{y_1}{y_2}}^{k}D_{
y}^1 \tilde{u}\|_{L^2( \tilde{\Omega})}^2 \mm{d}\tau\mbox{ for any }0\leq k\leq 3,
\end{aligned}$$
which, together with Lemma \ref{lem:0306}, implies \eqref{sumefork=1}.
\hfill$\Box$
\end{pf}
\subsection{Energy estimate of Gronwall-type}
Now, we are able to establish the energy estimate of Gronwall-type.
Putting  Lemmas \ref{lem:0305} and \ref{lem:0309} together, we get
\begin{equation*}\label{}\begin{aligned}
\int_0^t\left\{\|
 u\|_{H^4}^2+\left\| \frac{d}{dt}(\bar{e}\varrho+\bar{\rho}\theta)\right\|_{H^3}^2\right\}\mm{d}\tau
\lesssim \mathcal{E}_0^2 +\int_0^t\mathcal{R}\mm{d}\tau ,
\end{aligned}\end{equation*}
where the right-hand side can be bounded as follows, using Lemma \ref{lem:0304}.
$$\begin{aligned}\int_0^t\mathcal{R}\mm{d}\tau\lesssim  \mathcal{E}_0^2+ \int_0^t[ \|(\varrho,\theta)\|_{H^2}^2
+\mathcal{E}(\|u\|_{H^3}+\|u\|_{H^4}^2+\mathcal{E}^2)]\mm{d}\tau \end{aligned}.$$
Hence,
\begin{equation*}
\begin{aligned}
\int_0^t\left\{\|u\|_{H^4}^2+\left\| \frac{d}{dt}(\bar{e}\varrho+\bar{\rho}\theta)\right\|_{H^3}^2\right\}\mm{d}\tau
\lesssim \mathcal{E}_0^2+ \int_0^t\big[ \|(\varrho,\theta)\|_{H^2}^2
+\mathcal{E}(\|u\|_{H^3}+\mathcal{E}^2)\big] \mm{d}\tau  .
\end{aligned} \end{equation*}
Using interpolation inequality and Young's inequality, we further have
\begin{equation}    \label{asnoenet}
\begin{aligned}
\int_0^t\left\{\|u\|_{H^4}^2+\left\| \frac{d}{dt}(\bar{e}\varrho+\bar{\rho}\theta)\right\|_{H^3}^2\right\}\mm{d}\tau
\lesssim \mathcal{E}_0^2+ \int_0^t\big[ C_\epsilon\|(\varrho,u,\theta)\|_{L^2}^2
+\mathcal{E}^2(\epsilon+ \mathcal{E} )\big] \mm{d}\tau ,
\end{aligned} \end{equation}where the constant $C_\epsilon$ depends on $\epsilon$ and some physical parameters in \eqref{0105}.
On the other hand, by Lemmas \ref{lem:0301}--\ref{lem:0303}, \eqref{detetim1} and \eqref{momentum2}, we find that
$$\begin{aligned}
&\mathcal{E}^2(t)+\|(\varrho,\theta)_t\|_{H^2}^2+\|u_t(t)\|_{H^1}^2+\int_0^t\|u_{tt}\|^2_{L^2}
\mm{d}\tau\\
 & \lesssim \mathcal{E}^2_0+\int_0^t[\|u\|_{H^2}^2+
\mathcal{E}(\|u\|_{H^{4}} +\mathcal{E}^2)]\mm{d} \tau\\
 & \lesssim \mathcal{E}^2_0+\int_0^t[ C_\delta\| u \|_{L^2}^2+C_\delta\|u\|_{H^{4}}^2+
\mathcal{E}^2(\delta+\mathcal{E} )]\mm{d} \tau,
\end{aligned}$$where the constant $C_\delta$ depends on $\delta$ and some physical parameters in \eqref{0105}.
Consequently, in view of the above inequality and \eqref{asnoenet},  we obtain
\begin{equation}\label{ernygfor}\begin{aligned}
&\mathcal{E}^2(t)+\|(\varrho,\theta)_t\|_{H^2}^2+\|u_t(t)\|_{H^1}^2+\int_0^t\|u_{tt}\|^2_{L^2} \mm{d}\tau\\
& \leq C\left\{(1+C_\delta) \mathcal{E}_0^2 + \int_0^t\{C_\delta(1+C_\epsilon)\| (\varrho,u,\theta) \|_{L^2}^2+
\mathcal{E}^2[\delta +C_\delta \epsilon+(1+C_\delta)\mathcal{E} ]\}
\mm{d}\tau\right\}.
\end{aligned}
\end{equation}
In particular,   $\epsilon$, $\delta$ and $\mathcal{E}$ can be chosen so small that   $C[ \delta+C_\epsilon \epsilon  +(1+C_\delta  )\mathcal{E}]<\Lambda$ later on.

Now, let us recall that the local existence and uniqueness of solutions to the perturbed equations \eqref{0105} have been
established in \cite[Remark 6.1]{kawashima1984systems} for $\bar{\rho}$ and $\bar{e}$ being constants, while the global existence
and uniqueness of small solutions to the perturbed equations \eqref{0105} with heat conductivity have been shown
in \cite{MANTIJ321} for $(\bar{\rho},\bar{e})$ being close to a constant state.
By a slight modification in the proof of the local existence in \cite{kawashima1984systems} or \cite{MANTIJ321}, one can easily
obtain the existence and uniqueness of a local solution $(\varrho,v,\theta)\in C^0([0,T],H^3)$ to
the perturbed problem (\ref{0105})--(\ref{0107}) for some $T>0$. Moreover, this local solution satisfies
the above \emph{a priori} estimate (\ref{ernygfor}).
Therefore, we arrive at the following conclusion:
\begin{pro} \label{pro:0401} Assume that $(\bar{\rho},\bar{e})$ satisfies \eqref{0102} and \eqref{0104}.
For any given initial data $(\varrho_0, u_0,\theta_0)\in H^3$ satisfying the compatibility condition and
$$\inf_{x \in\Omega}\{\varrho_0 +\bar{\rho},\ \bar{\theta}_0+\bar{e}\}>0,$$
then there exist a $T^{\max}>0$ and a unique solution $(\varrho, u,\theta)\in C^0([0,T^{\max}),H^3)$ to the  Rayleigh-B\'enard problem  \eqref{0105}--\eqref{0107} satisfying $$\inf_{(0,T)\times\Omega}\{\varrho+\bar{\rho},\theta+\bar{e}\}>0,$$
where $T^{\max}$ denotes the maximal time of existence of the solution
$(\varrho, u,\theta)$.
Moreover, there is a sufficiently small constant ${\delta}^0_1\in (0,1]$, such that
if $\mathcal{E}(t)\leq {\delta}^0_1$ on some interval $[0,T]\subset [0,T^{\max})$, then the solution $(\varrho, u,\theta)$ satisfies
\begin{equation}\label{energyinequality}
\begin{aligned}
&\mathcal{E}^2(t)+\|(\varrho,\theta)_t(t)\|_{H^2}^2+\|u_t(t)\|_{H^1}^2+\int_0^t\|u_{tt}(\tau)\|^2_{L^2} \mm{d}\tau\\
&\quad \leq C\mathcal{E}_0^2+ \int_0^t(C\|(\varrho, u,\theta)(\tau)\|_{L^2}^2
+\Lambda\mathcal{E}^2(\tau)) \mm{d}\tau \;\;\mbox{ for any }t\in (0,T),
\end{aligned}\end{equation}
where the constant $C$ only depends on $ {\delta}_1^0$, $\Lambda$, $\Omega$ and  the known physical parameters in \eqref{0105}.
\end{pro}

\section{Nonlinear instability}\label{sec:04}
Now we are in a position to prove Theorem \ref{thm:0102}
by adopting and modifying the ideas in \cite{JFJSO2014,GYHCSDDC}. In view of Theorem \ref{thm:0101}, we can construct a (linear) solution
\begin{equation}\label{0501}
\left(\varrho^\mm{l},
{ u}^\mm{l},\theta^\mm{l}\right)=e^{{\Lambda t}}
\left(\bar{\varrho}_0,
\bar{ u}_0,\bar{\theta}_0\right)
\end{equation}
to the linearized problem \eqref{0106}--\eqref{0108} with the initial data $(\bar{\varrho}_0,\bar{ u}_0,\bar{\theta}_0)\in H^3 $.
Furthermore, this solution satisfies
\begin{equation}\label{wangweiwe16}
\|({\bar{u}}_{01},{\bar{u}}_{02})\|_{L^2}\|{\bar{u}}_{03}\|_{L^2}>0, \end{equation}
where $\bar{u}_{0i}$ stands for the $i$-th component of $\bar{ u}_0$ for $i=1$, $2$ and $3$.
In what follows,
$C_1,\cdots ,C_7$ will denote generic constants that may depend on $(\bar{\varrho}_0,
\bar{ u}_0,\bar{\theta}_0)$, $ {\delta}_1^0$, $\Lambda$, $\Omega$ and  the known physical parameters in \eqref{0105},
but are independent of $\delta$.

Obvious, we can not directly use the initial data of the linearized equations \eqref{0106}--\eqref{0108} as the one of the associated nonlinear
problem, since the linearized and nonlinear equations enjoy different compatibility conditions at the boundary.
To get around this obstacle,  we instead use the elliptic theory to construct initial data
of the nonlinear equations problem which are close to the linear growing modes.
\begin{lem}\label{lem:modfied} Let $(\bar{\varrho}_0,\bar{ u}_0, \bar{ {\theta}}_0)$ be the same as in \eqref{0501}.
 Then there exists a $ {\delta}^0_2\in (0,1)$ depending on $(\bar{\varrho}_0,\bar{ u}_0, \bar{ {\theta}}_0)$, such that
 for any $\delta\in (0, {\delta}^0_2)$, there is a $u_\mm{r}$ which may depend on $\delta$ and enjoys the following properties:
\begin{enumerate}
  \item[(1)] The modified initial data
  \begin{equation}\label{mmmode}(  {\varrho}_0^\delta,{ u}_0^\delta,{{\theta}}_0^\delta )
=\delta
   (\bar{\varrho}_0,\bar{ u}_0, \bar{ {\theta}}_0)
   +  \delta^2( \bar{\varrho}_0,{ u}_\mm{r}, \bar{ {\theta}}_0)
\end{equation}
satisfy ${ u}_0^\delta|_{\partial\Omega}= 0$ and the compatibility condition:
$$\big\{(\varrho_0^\delta+\bar{\rho}){ u}_0^\delta\cdot\nabla {u}_0^\delta+a\nabla [{({\varrho}_0^\delta
+\bar{\rho})(\theta_0^\delta+\bar{e})} -{\bar{\rho}\bar{e}}]-\mu\Delta{ u}_0^\delta
-\mu_0\nabla\mm{div} u_0^\delta+g\varrho_0^\delta  e_3\big\}|_{\partial\Omega}=0.  $$
  \item[(2)] $ { u}_\mm{r} $ satisfies the following estimate:
$$ \| { u}_\mm{r}\|_{H^{3}} \leq C_1, $$
where the constant $C_1$ depends on $\|(\bar{\varrho}_0,\bar{ u}_0, \bar{ {\theta}}_0)\|_{H^3}$
and other physical parameters, but is independent of $\delta$.
\end{enumerate}
\end{lem}
\begin{pf}
Notice that $(\bar{\varrho}_0,\bar{ u}_0, \bar{ {\theta}}_0)$ satisfies
 $$\bar{ u}_0|_{\partial\Omega}= 0,\quad [a\nabla(\bar{e} \bar{\varrho}_0
  +\bar{\rho}\bar{\theta}_0)-\mu\Delta\bar{ u}_0-\mu_0\nabla \mm{div}\bar{ u}_0+g\bar{\varrho}_0 e_3]|_{\partial\Omega}=0.$$
 Hence, if the modified initial data satisfy \eqref{mmmode}, then we expect ${u}_\mm{r}$ to satisfy the following problem:
\begin{equation} \left\{\begin{array}{l}
 \mu\Delta{ u}_\mm{r} +\mu_0\nabla\mm{div} u_\mm{r} -\delta^2\varrho_0^{**}{ u}_\mm{r}  \cdot\nabla
{u}_\mm{r}-\delta\varrho_0^{**}(\bar{ u}_0\cdot\nabla {u}_\mm{r}+{ u}_\mm{r}\cdot\nabla\bar{ u}_0)  \\
\quad=a\nabla(\bar{e} \bar{\varrho}_0   +\bar{\rho}\bar{\theta}_0)+g\bar{\varrho}_0 e_3+
\varrho_0^{**}\bar{ u}_0\cdot\nabla \bar{ u}_0-a\nabla ({\varrho}_0^*{\theta}_0^*):=F(\bar{\varrho}_0,\bar{u}_0,\bar{\theta}_0), \\[1mm]
  u_\mm{r}|_{\partial\Omega} = 0
 \end{array} \right.   \label{js2} \end{equation}
where $\varrho_0^{*}:=(1 +\delta )\bar{\varrho}_0$, $\theta_0^{*}=(1 +\delta )\bar{\theta}_0$ and
$\varrho_0^{**}:=(\varrho_0^\delta+\bar{\rho})=(\delta +\delta^2 )\bar{\varrho}_0+\bar{\rho}$.
Thus the modified initial data naturally satisfy the compatibility condition.

Next we shall look for a solution $u_\mm{r}$ to the boundary problem (\ref{js2}) when $\delta$ is sufficiently small.
We begin with the linearization of (\ref{js2}) which reads as
\begin{equation}\label{elliequation}
  \begin{aligned}&\mu\Delta{ u}_\mm{r} +\mu_0\nabla\mm{div} u_\mm{r}
  =F(\bar{\varrho}_0,\bar{ u}_0,\bar{\theta}_0)+\delta^2\varrho_0^{**}{v} \cdot\nabla
{v}+\delta\varrho_0^{**}(\bar{ u}_0\cdot\nabla {v}+{v}\cdot\nabla\bar{ u}_0)  \end{aligned}
\end{equation}
with boundary condition
\begin{equation}\label{boundery0}
{ u}_\mm{r}|_{\Omega}= 0.
\end{equation}

Let $v\in H^{3}$, then it follows from the elliptic theory that there is a solution ${u}_\mm{r}$ of
\eqref{elliequation}--\eqref{boundery0} satisfying
$$ \begin{aligned}
\|{ u}_\mm{r}\|_{H^{3}} \leq & \|F(\bar{\varrho}_0,\bar{u}_0,\bar{\theta}_0) +\delta^2\varrho_0^{**}{v}\cdot\nabla
{v}+\delta\varrho_0^{**}(\bar{u}_0\cdot\nabla {v}+{v}\cdot\nabla\bar{u}_0)\|_{H^1}\\[1mm]
\leq & C_{\mm{m}}(1+\|(\bar{\varrho}_0, \bar{u}_0, \bar{\theta}_0)\|_{H^{2}}^2+\delta^2 \|v\|_{H^{2}}^2).
\end{aligned}  $$
Now, we take $C_1=C_{\mm{m}}(2+\|(\bar{\varrho}_0,\bar{ u}_0,\bar{\theta}_0)\|_{H^{2}}^2)$ and $\delta\leq \min\{C^{-1}_1,1\}$.
Then for any $\| v\|_{H^{3}}^2\leq C_1$, one has
$$\begin{aligned}\|{ u}_\mm{r}\|_{H^{3}}
\leq C_1.
\end{aligned}
$$
Therefore we can construct an approximate function sequence ${u}_\mm{r}^n$, such that
\begin{equation*}  \begin{aligned}
& \mu\Delta{ u}_\mm{r}^{n+1} +\mu_0\nabla\mm{div} u_\mm{r}^{n+1}
-\delta^2\varrho_0^{**}{ u}^n_\mm{r} \cdot\nabla {u}^n_\mm{r}-\delta\varrho_0^{**}(\bar{ u}_0\cdot\nabla
{u}^n_\mm{r} +{u}^n_\mm{r}\cdot\nabla\bar{ u}_0)=F(\bar{\varrho}_0, \bar{u}_0, \bar{\theta}_0), \end{aligned}
\end{equation*}
and for any $n$, $$\|{ u}_\mm{r}^n\|_{H^{3}}
\leq C_1,\quad \|{ u}_\mm{r}^{n+1}-u_\mm{r}^{n}\|_{H^{3}}
\leq C_2\delta\|{ u}_\mm{r}^{n}-u_\mm{r}^{n-1}\|_{H^{3}}$$
for some constant $C_2$ independent of $\delta$ and $n$.
Finally, we choose a $\delta$ sufficiently small so that $C_2\delta<1$, and then use a compactness argument
to get a limit function which solves the nonlinear boundary problem (\ref{js2}). Moreover $\|{u}_\mm{r}\|_{H^{3}} \leq C_1$.
Thus we have proved Lemma \ref{lem:modfied}.    \hfill$\Box$
\end{pf}

Let $({\varrho}_0^\delta,{ u}_0^\delta,{{\theta}}_0^\delta )$ be constructed as in Lemma \ref{lem:modfied}. Then there is
a constant $$C_3\geq \max\{1, \|\left(\bar{\varrho}_0,
\bar{ u}_0,\bar{\theta}_0\right)\|_{L^2}\}$$ depending on $(\bar{\varrho}_0,\bar{u}_0,\bar{{\theta}}_0)$,
such that for any $\delta\in (0,{\delta}_2^0)\subset (0,1)$,
$$\mathcal{E}({\varrho}_0^\delta,{ u}_0^\delta,{{\theta}}_0^\delta )\leq C_3\delta,  $$
where $\mathcal{E}$ is defined by \eqref{enerdienf}.
Recalling $\inf_{ x\in\Omega}\{\bar{\rho},\bar{e}\}>0$ and the embedding theorem $H^2\hookrightarrow L^\infty$,
we can choose a sufficiently small $\delta$, such that
\begin{equation}\label{inferfds}
\inf_{ x\in\Omega}\{\varrho_0^\delta+\bar{\rho}, \theta_0^\delta+\bar{e}\}>0.\end{equation}
Hence, by virtue of Proposition \ref{pro:0401}, there is a $ {\delta}^0_3\in (0, {\delta}^0_2)$,
such that for any $\delta<{\delta}^0_3$, there exists a unique local solution
$(\varrho^\delta, u^\delta, \theta^\delta)\in C^0([0,T],H^3)$ to \eqref{0105} and \eqref{0107}, emanating
from the initial data $(\varrho_0^\delta, u_0^\delta,\theta_0^\delta)$. Moreover,
\eqref{inferfds} holds for any $\delta$ satisfying $\mathcal{E}({\varrho}_0^\delta,{u}_0^\delta,{{\theta}}_0^\delta )\leq C_3{\delta}^0_3$.
Let $C>0$ and ${\delta}^0_1>0$ be the same constants as in Proposition \ref{pro:0401} and
$\delta_0=\min\{{\delta}^0_3, {\delta}^0_1/C_3\}$. Let $\delta\in (0,\delta_0)$ and
 \begin{equation}\label{times}
 T^{\delta}=\frac{1}{\Lambda}\mm{ln}\frac{2\varepsilon_0}{\delta}>0,\;\quad\mbox{i.e.,}\;\;
 \delta e^{\Lambda T^\delta}=2\varepsilon_0, \end{equation}
where $\varepsilon_0\leq 1$, independent of $\delta$, is sufficiently small and will be fixed later.
In what follows, we denote $\mathcal{E}_\delta(t):={\mathcal{E}}(\varrho^\delta,{ u}^\delta,{\theta}^\delta )(t)$.

Define
 \begin{equation*}
T^*=\sup\left\{t\in (0,T^{\max})\left|~{\mathcal{E}}_\delta(t)\leq C_3{\delta_0}\right.\right\}>0\end{equation*}
 and
 \begin{equation*}    T^{**}=\sup\left\{t\in (0,T^{\max})\left|~\left\|\left(\varrho^\delta,
{u}^\delta,\theta^\delta\right)(t)\right\|_{{L}^2}\leq 2\delta C_3e^{\Lambda t}\right\}>0\right., \end{equation*}
where $T^{\mm{max}}$ denotes the maximal time of existence of the solution $(\varrho^\delta,{u}^\delta,\theta^\delta)$.
Obviously, $T^*$ and $T^{**}$ may be finite, and furthermore,
 \begin{eqnarray}\label{0502n1}
&&\mathcal{E}_\delta(T^*)=C_3{\delta_0}\quad\mbox{ if }T^*<\infty ,  \\[1mm]
\label{0502n111}  &&  \left\|\left(\varrho^\delta, { u}^\delta,\theta^\delta\right)(T^{**})\right\|_{{L}^2}
=2\delta C_3e^{\Lambda T^{**}}\quad\mbox{ if }T^{**}<T^{\max}.
\end{eqnarray}
Then for all $t\leq \min\{T^\delta,T^*,T^{**}\}$, we deduce from the estimate \eqref{energyinequality} and the
definition of $T^*$ and $T^{**}$ that
 \begin{equation*}\begin{aligned}
&{\mathcal{E}}^2_\delta(t) +\|(\varrho^\delta,\theta^\delta)_t(t)\|_{H^2}^2 +\|u_t^\delta(t)\|_{H^1}^2
+\int_0^t\|u_{tt}^\delta\|^2_{L^2} \mm{d}\tau  \\
& \leq C[ \mathcal{E}^2({\varrho}_0^\delta,{ u}_0^\delta,{{\theta}}_0^\delta ) +2 C_3^2\delta^2e^{2\Lambda t}/\Lambda]
+\Lambda\int_0^t{\mathcal{E}}^2_\delta(\tau) \mm{d}\tau  \\
& \leq  C_4\delta^2e^{2\Lambda t} +\Lambda\int_0^t{\mathcal{E}}^2_\delta(\tau) \mm{d}\tau
   \end{aligned}   \end{equation*}
for some constant $C_4>0$. Thus, applying Gronwall's inequality, one concludes
 \begin{equation}\begin{aligned}\label{0503}
{\mathcal{E}}^2_\delta(t) +\|(\varrho^\delta,\theta^\delta)_t(t)\|_{H^2}^2
+ \|u_t^\delta(t)\|_{H^1}^2 +\int_0^t\|u_{tt}^\delta\|^2_{L^2} \mm{d}\tau
\leq C_5\delta^2e^{2\Lambda t}       \end{aligned}   \end{equation}
for some constant $C_5>0$.

Let $(\varrho^{\mathrm{d}}, {u}^{\mathrm{d}}, {\theta}^{\mathrm{d}})=(\varrho^{\delta},
{u}^{\delta},{\theta}^{\delta})-\delta(\varrho^{\mathrm{l}},{u}^{\mathrm{l}},{\theta}^{\mathrm{l}})$.
Noting that $(\varrho^\mm{a}_\delta, u^{\mm{a}}_\delta, \theta^{\mm{a}}_\delta):= \delta(\varrho^{\mm{l}}, u^{\mm{l}}, \theta^{\mm{l}})$
is also a solution to the linearized problem \eqref{0106}--\eqref{0108} with the initial data
$\delta(\bar{\varrho}_0, \bar{u}_0, \bar{\theta}_0)\in H^3$, we find that $(\varrho^{\mathrm{d}}, {u}^{\mathrm{d}}, {\theta}^{\mathrm{d}})$
satisfies the following error equations:
\begin{equation}\label{h0407}\left\{\begin{array}{ll}
  \varrho_t^{\mathrm{d}} + \mm{div}(\bar{\rho} u^{\mathrm{d}})= -\mm{div}(\varrho^{\delta} u^{\delta})
  :=N^{\varrho}(\varrho^\delta, u^{\delta}):=N^\varrho_\delta, \\[2mm]
  \bar{\rho} u_t^{\mathrm{d}} + a\nabla(\bar{e} {\varrho}^{\mathrm{d}}
+\bar{\rho}\theta^{\mathrm{d}})-\mu\nabla\mm{div} u^{\mathrm{d}}-\mu_0\Delta u^{\mathrm{d}}+g\varrho^{\mathrm{d}} e_3 \\
\qquad =   -( \varrho^{\delta}+\bar{\rho}) u^{\delta}\cdot\nabla u^{\delta} -\varrho^\delta u^{\delta}_t
-a\nabla (\varrho^\delta \theta^{\delta}):=N^u(\varrho^\delta,u^{\delta},\theta^{\delta}):=N^u_\delta, \\[2mm]
 \theta_t^{\mathrm{d}}+\bar{e}'u^{\mathrm{d}}_3 +a\bar{e}\mm{div}u^{\mathrm{d}}=
[{\mu}|\nabla u^{\delta}+\nabla (u^{\delta})^\mm{T}|^2/2 + \lambda(\mm{div}u^{\delta})^2]/(\varrho^{\delta}+\bar{\rho}) \\[2mm]
\qquad\quad\qquad\quad\qquad \qquad  -u^{\delta}\cdot\nabla \theta^{\delta}
-a\theta^{\delta}\mm{div}u^{\delta}:=N^e(\varrho^\delta,u^{\delta}, \theta^{\delta}):=N^\theta_\delta,
 \end{array}\right.\end{equation}
 with initial data $(\varrho^{\mathrm{d}}(0),{u}^{\mathrm{d}}(0),\theta^{\mathrm{d}}(0))=
 \delta^2(\bar{\varrho}_0,{u}_\mm{r},\bar{{\theta}}_0)$ and boundary condition $u^{\mm{d}}|_{\partial\Omega}=0$.

  Next, we shall establish the error estimate for $(\varrho^{\mm{d}},u^{\mm{d}},\theta^{\mm{d}})$ in $L^2$-norm.
   \begin{lem}\label{erroestimate}  There is a constant $C_6$, such that for all $t\leq \min\{T^\delta,T^*,T^{**}\}$,
    \begin{equation}\label{ereroe}
\begin{aligned}
  \| (\varrho^{\mm{d}},u^{\mm{d}},\theta^{\mm{d}})(t)\|^2_{L^2} \leq C_6\delta^3\theta^{3\Lambda t}.
 \end{aligned}  \end{equation}
   \end{lem}
\begin{pf}
We differentiate the linearized  momentum equations \eqref{h0407}$_{2}$ in time, multiply the resulting
equations by $u_t$ in $L^2(\Omega)$, and use the equations \eqref{h0407}$_{1}$ and \eqref{h0407}$_{3}$ to deduce
\begin{equation}\label{nnn0314P}
\begin{aligned}
&\frac{\mm{d}}{\mm{d}t}
\int \left\{\bar{\rho}| u_t^{\mathrm{d}}|^2
-g\bar{\rho}'|{u}_3^{\mathrm{d}}|^2+
[(1+a)\bar{p}\mm{div} u^{\mathrm{d}}
-2g\bar{\rho}{u}_3^{\mathrm{d}}]\mm{div}{ u}^{\mathrm{d}}\right\}\mathrm{d} x\\
&=-
2\mu\int |\nabla u^{\mathrm{d}}_t|^2\mm{d} x-
2\mu_0\int |\mm{div} u_t^{\mathrm{d}}|^2
\mm{d} x\\
&\quad+ 2\int [\partial_tN^u_\delta-
gN^\varrho_\delta  e_3-  a\nabla(\bar{e}N^\varrho_\delta +\bar{\rho} N^\theta_\delta)]\cdot
 u_t^{\mathrm{d}}\mm{d} x
.
\end{aligned}\end{equation}
Thanks to \eqref{sharprate}, one has
\begin{equation*}\label{0302}\begin{aligned}
&\int \{g\bar{\rho}'|{u}^{\mathrm{d}}_3|^2 +[2g\bar{\rho}u_3^{\mathrm{d}}-(1+a)\bar{p}
\mm{div}{ u}^{\mathrm{d}}]\mm{div}{ u}^{\mathrm{d}}\} \mm{d} x\\
& \quad \leq\Lambda\int\left(\mu|\nabla {u}^{\mathrm{d}}|^2+\mu_0|\mm{div}
{u}^{\mathrm{d}}|^2\right)\mm{d}x +\Lambda^2{\int\bar{\rho}|{ u}^{\mathrm{d}}|^2\mathrm{d} x}.
\end{aligned}\end{equation*}
Thus, integrating (\ref{nnn0314P}) in time from $0$ to $t$, we get
\begin{equation}\label{0314}
\begin{aligned}
&\|\sqrt{\bar{\rho}} u_t^\mm{d}(t)\|^2_{L^2}+2\int_0^t(\mu\|\nabla   u_\tau ^\mm{d}\|^2_{L^2}
+\mu_0\|\mm{div}  u_\tau^\mm{d}\|^2_{L^2}) \mm{d}\tau  \\
& \leq I_1^0 + {\Lambda^2} \|\sqrt{\bar{\rho}} u^\mm{d}(t)\|_{L^2}
+ {\Lambda}\mu\|\nabla u^\mm{d}(t)\|^2_{L^2}+ {\Lambda}\mu_0\|\mm{div} u^\mm{d}(t)\|^2_{L^2}\\
& \quad + 2\int_0^t\int [\partial_\tau N^u_\delta - gN^\varrho_\delta  e_3
-  a\nabla(\bar{e}N^\varrho_\delta +\bar{\rho} N^\theta_\delta)]\cdot u^{\mathrm{d}}_\tau\mm{d} x\mm{d}\tau,
\end{aligned}\end{equation}
where $$I_1^0=\left\{\int \left\{\bar{\rho}| u_t^{\mathrm{d}}|^2
-g\bar{\rho}'({u}_3^{\mathrm{d}})^2+
[(1+a)\bar{p}\mm{div} u^{\mathrm{d}}
-2g\bar{\rho}{u}_3^{\mathrm{d}}]\mm{div}{ u}^{\mathrm{d}}\right\}\mathrm{d} x\right\}\bigg|_{t=0}.$$

Using Newton-Leibniz's formula and Cauchy-Schwarz's inequality, we find that
 \begin{equation}\begin{aligned}\label{0316}
& \Lambda (\mu\|\nabla u^\mm{d}(t)\|_{L^2}^2+\mu_0\|\mm{div} u^\mm{d}(t)\|^2_{L^2})  \\
& =I_2^0+ 2\Lambda\int_0^t\int_{\Omega}\left(\mu\sum_{1\leq i,j\leq 3}\partial_{x_i} u_{j\tau}^\mm{d}\partial_{x_i} u_{j\tau}^\mm{d} \mm{d} x\mathrm{d}\tau
+\mu_0\mm{div} u_\tau^\mm{d}\mm{div} u^\mm{d}\right)\mm{d} x\mathrm{d}\tau \\
& \leq I_2^0+\int_0^t(\mu\|\nabla u_\tau^\mm{d}\|_{L^2}^2 +\mu_0\|\mm{div} u_\tau^\mm{d}\|^2_{L^2})
\mathrm{d}\tau +\Lambda^2\int_0^t(\mu\|\nabla u^\mm{d}\|_{L^2}^2+\mu_0\|\mm{div} u^\mm{d}
\|^2_{L^2})\mathrm{d}\tau,       \end{aligned}\end{equation}
where $I_2^0=\Lambda (\mu\|\nabla u^\mm{d}(0)\|_{L^2}^2+\mu_0\|\mm{div} u^\mm{d}(0)\|^2_{L^2})$ and $u_{j\tau}^\mm{d}$ denotes the $j$-th component of
 $u_{\tau}^\mm{d}$ .
On the other hand,
\begin{equation}\begin{aligned}\label{0317}
\Lambda\frac{\mm{d}}{\mm{d}t}\|\sqrt{\bar{\rho}} u^\mm{d}(t)\|^2_{L^2}=2\Lambda\int
\bar{\rho} u^\mm{d}(t)\cdot  u^\mm{d}_t(t)\mm{d} x\leq\|\sqrt{\bar{\rho}}  u_t^\mm{d}(t)\|^2_{L^2}
+\Lambda^2\|\sqrt{\bar{\rho}} u^\mm{d}(t)\|^2_{L^2}.
\end{aligned}\end{equation}
Hence, putting (\ref{0314})--(\ref{0317}) together, we obtain the differential inequality
\begin{equation}\label{12safd}
\begin{aligned}
& \frac{\mm{d}}{\mm{d}t}\|\sqrt{\bar{\rho}} u^\mm{d}(t)\|^2_{L^2}+ \mu\| \nabla u^\mm{d}(t)\|_{L^2}^2+
\mu_0\| \mm{div} u^\mm{d}(t)\|^2_{L^2}  \\
&\leq2\Lambda\left[ \|\sqrt{\bar{\rho}} u^\mm{d}\|^2_{L^2} +\int_0^t(\mu\| \nabla u^\mm{d}\|_{L^2}^2 +
{\mu}_0\|\mm{div} u^\mm{d}\|_{L^2}^2) \mathrm{d}s\right]  \\
& \quad +\frac{I_1^0+2I^0_2}{\Lambda}+ \frac{2}{\Lambda}\int_0^t \int [\partial_\tau N^u_\delta-
gN^\varrho_\delta  e_3-  a\nabla(\bar{e}N^\varrho_\delta +\bar{\rho} N^\theta_\delta)]\cdot
u^{\mathrm{d}}_\tau\mm{d} x\mm{d}\tau.
\end{aligned}
\end{equation}

Next, we control the last two terms on the right hand of \eqref{12safd}.
Noting that \begin{equation}\label{timesd}
\delta e^{\Lambda t}\leq 2\varepsilon_0\leq 2 \quad\mbox{for any }t\leq \min\{T^\delta,T^*,T^{**}\},
\end{equation}
we utilize \eqref{0503} and \eqref{0501}, H\"older's inequality and Sobolev's embedding theorem to infer that
  \begin{equation}\label{esitmaeintial1}
\begin{aligned}
& \left|2\int_0^t \int [\partial_\tau N^u_\delta - gN^\varrho_\delta e_3
- a\nabla(\bar{e}N^\varrho_\delta +\bar{\rho} N^\theta_\delta)]\cdot
  u^{\mathrm{d}}_\tau\mm{d} x\mm{d}\tau\right|        \\
& \lesssim \int_0^t (\|(N^\varrho_\delta,N_\delta^\theta)\|_{H^1}+\|\partial_\tau N^u_\delta\|_{L^2})(\|u_\tau^\mm{a}
\|_{L^2}+\|u_\tau^\delta\|_{L^2})\mm{d}\tau   \\
& \lesssim  \int_0^t (\delta^3e^{3\Lambda \tau}+\delta^2e^{2\Lambda \tau}
+\delta e^{\Lambda \tau}\|u_{\tau\tau}^\delta\|_{L^2})\delta e^{\Lambda \tau}\mm{d}\tau \\
& \lesssim  \delta^3e^{3\Lambda t}+ \delta^4e^{4\Lambda t}\lesssim   \delta^3e^{3\Lambda t},
 \end{aligned}  \end{equation}
 and
 \begin{equation}\label{esitmaeintial2}
\begin{aligned}
({I_1^0 +2I^0_2})/{\Lambda}\lesssim&(\|\sqrt{\bar{\rho}} u_t^\mm{d}\|^2_{L^2}+\|\nabla u^\mm{d}_0\|_{L^2}^2+\|u^\mm{d}_3\|_{L^2}^2)|_{t=0}\\
\lesssim& [\|(\varrho^\mm{d},\theta^\mm{d})\|_{H^1}^2+\|u^\mm{d}\|_{H^2}^2+\mathcal{E}_\delta^2
(\mathcal{E}_\delta^2+\mathcal{E}_\delta^4+\|u_t^\delta\|_{L^2}^2)]|_{t=0}
\\
\lesssim & \delta^4(\|(\bar{\varrho}_\mm{0},\bar{\theta}_{0})\|_{H^1}^2
+\|u_{\mm{r}}\|_{H^2}^2)+\delta^2 e^{2\Lambda t}(\delta^2 e^{2\Lambda t}+\delta^4 e^{4\Lambda t})\lesssim \delta^3 e^{3\Lambda t} .
 \end{aligned}  \end{equation}
Thus, substituting (\ref{esitmaeintial2}) and (\ref{esitmaeintial1}) into (\ref{12safd}), we obtain
$$  \begin{aligned}
& \partial_t \|\sqrt{\bar{\rho}} u^{\mathrm{d}}(t)\|^2_{L^2}+
\mu\|\nabla u^{\mathrm{d}}(t)\|_{L^2}^2 + {\mu_0}\|\mm{div} u(t)\|^2_{L^2}  \\
& \leq 2\Lambda\left[ \|\sqrt{\bar{\rho}} u^{\mathrm{d}}(t)\|^2_{L^2} +\int_0^t(\mu\|
 \nabla u^\mathrm{d}\|_{L^2}^2 + {\mu}_0\|\mm{div} u^{\mathrm{d}}\|^2_{L^2})
\mathrm{d}\tau\right]+C_7\delta^3e^{3\Lambda t}.
\end{aligned}$$
Applying Gronwall's inequality to the above inequality, one obtains
 \begin{equation}\label{estimerrvelcoity}
\begin{aligned}
  \|\sqrt{\bar{\rho}} u^{\mathrm{d}}(t)\|^2_{L^2}+
\int_0^t({\mu}\|\nabla u^{\mathrm{d}}\|^2_{L^2}+{\mu}_0\|\mm{div} u^{\mathrm{d}}\|^2_{L^2})\mm{d}\tau
\lesssim \delta^3e^{3\Lambda t}+\delta^4\|\sqrt{\bar{\rho}}u_{\mm{r}}\|_{L^2}^2\lesssim \delta^3e^{3\Lambda t}
 \end{aligned}  \end{equation}
 for all $t\leq \min\{T^\delta,T^*,T^{**}\}$.
 Thus, making use of \eqref{0314}, \eqref{0316} and \eqref{esitmaeintial1}--\eqref{estimerrvelcoity}, we deduce that
\begin{equation}\label{inequalemee}\begin{aligned}
&\frac{1}{\Lambda}\|\sqrt{\bar{\rho}}  u_t^{\mathrm{d}}(t)\|^2_{L^2}+
{\mu}\|\nabla u^{\mathrm{d}}(t)\|_{L^2}^2 +\mu_0\|\mm{div} u^{\mathrm{d}} (t)\|^2_{L^2}\\
& \leq   {\Lambda}\|\sqrt{\bar{\rho}} u^{\mathrm{d}}(t)\|^2_{L^2}+2 {\Lambda}\int_0^t({\mu}\|
\nabla u^{\mathrm{d}}\|_{L^2}^2+ {\mu}_0\|\mm{div} u^{\mathrm{d}}\|^2_{L^2})\mm{d}\tau\\
&\quad +\frac{I_1^0+2I^0_2}{\Lambda}+ \frac{2}{\Lambda}\int_0^t \int [\partial_\tau N^u_\delta-
gN^\varrho_\delta  e_3-  a\nabla(\bar{e}N^\varrho_\delta +\bar{\rho} N^\theta_\delta)]\cdot
u^{\mathrm{d}}_\tau\mm{d} x\mm{d}\tau\lesssim \delta^3e^{3\Lambda t}.
\end{aligned}\end{equation}
which, together with Poincar\'e's inequality and the estimates \eqref{estimerrvelcoity}, yields
\begin{eqnarray}\label{uestimate1n}
\| u^{\mathrm{d}}(t)\|_{H^1 }^2+\| u_t^{\mathrm{d}}(t)\|^2_{L^2 }+
\int_0^t\|\nabla u^{\mathrm{d}}\|^2_{L^2}\mm{d}\tau \lesssim \delta^3e^{3\Lambda t}.
\end{eqnarray}

Finally, using the equations \eqref{h0407}$_1$ and \eqref{h0407}$_2$, and the estimates \eqref{timesd} and \eqref{uestimate1n},
we find that
\begin{equation*}\begin{aligned}\label{}
\|(\varrho^{\mathrm{d}},\theta^{\mathrm{d}})(t)\|_{L^2}\leq & \delta^2\|(\bar{\varrho}_0,\bar{\theta}_0)\|_{L^2}+\int_0^t
\|(\varrho^{\mathrm{d}},\theta^{\mathrm{d}})_\tau\|_{L^2}\mm{d}\tau \\
\lesssim &\delta^{2}+\int_0^t(\| u^{\mm{d}}\|_{H^1}+\|(N^\varrho_\delta,N^\theta_\delta)\|_{L^2})\mm{d}\tau \\
\lesssim & \delta^2+\int_0^t(\delta^\frac{3}{2}e^{\frac{3\Lambda}{2}\tau}+\mathcal{E}_\delta^2(\tau))\mm{d}\tau\lesssim \delta^\frac{3}{2}e^{\frac{3\Lambda}{2} t}.
\end{aligned}\end{equation*}
Putting the previous estimates together, we get \eqref{ereroe} immediately.
This completes the proof of Lemma \ref{erroestimate}.  \hfill$\Box$
 \end{pf}

 Now, we claim that
\begin{equation}\label{n0508}
T^\delta=\min\left\{T^\delta,T^*,T^{**}\right\},
 \end{equation}
provided that small $\varepsilon_0$ is taken to be
 \begin{equation}\label{defined}
\varepsilon_0=\min\left\{\frac{{C_3\delta_0}}
{4\sqrt{C_5}},\frac{C_3^2}{8C_6},\frac{m_0^2}{2^3C_6},1 \right\}>0,
 \end{equation}
 where we have denoted $m_0=\min\{\|(\bar{u}_{01},\bar{u}_{02})\|_{L^2},\|\bar{u}_{03}\|_{L^2}\}>0$ due to \eqref{wangweiwe16}.

 Indeed, if $T^*=\min\{T^{\delta},T^*,T^{**}\}$, then $T^*<\infty$. Moreover, from \eqref{0503} and \eqref{times} we get
 \begin{equation*}
{\mathcal{E}}_\delta(T^*)\leq \sqrt{C_5}\delta e^{\Lambda T^*}
\leq \sqrt{C_5}\delta e^{\Lambda T^\delta}=2\sqrt{C_5}\varepsilon_0<C_3{\delta_0},
 \end{equation*}
 which contradicts with \eqref{0502n1}. On the other hand, if $T^{**}=\min\{T^{\delta},T^*,T^{**}\}$, then $T^{**}<T^*\leq T^{\mm{max}}$.
Moreover, in view of \eqref{0501}, \eqref{times} and \eqref{ereroe}, we see that
 \begin{equation*}\begin{aligned}
 \left\|\left(\varrho^\delta, { {u}}^\delta,{ {\theta}}^\delta \right)(T^{**})\right\|_{L^2}
\leq  & \left\|\left(\varrho^\mm{a}_{\delta}, { {u}}^\mm{a}_{\delta},{ {\theta}}_\delta^{\mm{a}}
\right)(T^{**})\right\|_{L^2} +\left\|\left(\varrho^{\mathrm{d}},
{ {u}}^{\mathrm{d}},{ {\theta}}_\delta^{\mm{d}} \right)(T^{**})\right\|_{L^2} \\
\leq  &\delta \left\|\left(\varrho^\mm{l}, { {u}}^{\mm{l}},{ {\theta}}^{\mm{l}}
\right)(T^{**})\right\|_{L^2}+\sqrt{C_6}\delta^{3/2}e^{3\Lambda T^{**}/2} \\
\leq & \delta C_3e^{\Lambda T^{**}}+\sqrt{C_6}\delta^{3/2} e^{3\Lambda T^{**}/2} \\
\leq & \delta e^{\Lambda T^{**}}(C_3+\sqrt{2C_6\varepsilon_0})
<2\delta C_3  e^{\Lambda T^{**}},
 \end{aligned} \end{equation*}
which also contradicts with \eqref{0502n111}. Therefore, \eqref{n0508} holds.

 Since $T^\delta=\min\left\{T^\delta,T^*,T^{**}\right\}$, \eqref{ereroe} holds for $t= T^\delta$,
thus we again use \eqref{defined} and \eqref{ereroe} with $t= T^\delta$ to deduce that
 \begin{equation*}\begin{aligned}
 \|u_3^{\delta}(T^\delta)\|_{L^2}\geq &
\|u^{\mathrm{a}}_{3\delta}(T^{\delta})\|_{L^2}-\|u_3^{\mm{d}}(T^{\delta})\|_{L^2}
= \delta e^{\Lambda T^\delta}\|\bar{u}_{03}\|_{L^2}-\|u_3^{\mm{d}}(T^{\delta})\|_{L^2} \\
 \geq & \delta e^{\Lambda T^\delta}\|\bar{u}_{03}\|_{L^2}-\sqrt{C_6}\delta^{3/2}e^{3\Lambda^* T^{\delta}/2}
 \geq  2m_0\varepsilon_0-\sqrt{C_6}(2\varepsilon_0)^{3/2}  \geq m_0\varepsilon_0,
 \end{aligned}      \end{equation*}
 where $u^{\delta}_{3}(T^{\delta})$  denote the third component of $u^{\delta}(T^{\delta})$.
 Similar, we also have
 $$\|(u_1^{\delta},u_2^{\delta})(T^\delta)\|_{L^2}\geq m_0\varepsilon_0.$$
This completes the proof of Theorem \ref{thm:0102} by defining $\varepsilon=m_0\varepsilon_0$.
In addition, if
$\bar{\rho}'\geq 0$, then the function $\bar{\rho}_0$ constructed in \eqref{0501} satisfies $\|\bar{\rho}_0\|_{L^2}>0$.
Thus we also obtain $\|\varrho^\delta(T^\delta)\|_{L^2}\geq m_0\varepsilon_0$,
if we define $m_0=\min\{\|\bar{\varrho}_{0}\|_{L^2}, \|(\bar{u}_{01},
\bar{u}_{02})\|_{L^2},\|\bar{u}_{03}\|_{L^2}\}>0$. Hence, the assertion in Remark \ref{strongconden} holds.
\setcounter{equation}{0}  
\section*{Appendix}\label{sec:03}
\appendix
  \renewcommand{\theequation}{A.\arabic{equation}}
  In this section we show that $\Lambda$ defined by \eqref{sharprate} is the sharp growth rate for any
solutions to the linearized problem \eqref{0106}--\eqref{0108}. Since the density varies for a compressible fluid,
the spectrums of the linearized solution operator are difficult to analyze in comparison with an incompressible fluid,
and it is hard to obtain the largest growth rate of the solution operator in some Sobolev space in the usual way.
Here we exploit energy estimates as in \cite{HHJGY,JFJSO2014} to show that $e^{\Lambda t}$
is indeed the sharp growth rate for $(\varrho,u,\theta)$ in $H^2$-norm.
\begin{pro} \label{grwothe}  Assume that the assumption of Theorem \ref{thm:0101} is satisfied. Let
$(\varrho, u,\theta)$ solve the  linearized Rayleigh-B\'enard problem \eqref{0106}--\eqref{0108}.
Then, we have the following estimates.
\begin{eqnarray}
&&   \hspace{-15mm}
\|(\varrho,\theta)(t)\|_{L^2}^2+ \| u(t)\|_{H^1 }^2+\| u_t(t)\|^2_{L^2 }+
\int_0^t\|\nabla u(s)\|^2_{L^2}\mm{d}s \leq Ce^{2\Lambda t}(\|(\varrho_0,\theta_0)\|_{L^2}^2+\| u_0\|_{H^2}^2), \label{uestimate}  \\
&&  \hspace{-15mm}  \|(\varrho,u,\theta)(t)\|_{H^2}^2+\int_0^t\left[\|
 u\|_{H^3}^2+\left\| (\bar{e}\rho_t+\bar{\rho}\theta_t)\right\|_{H^2}^2\right]\mm{d}\tau
\leq Ce^{2\Lambda t}\|(\varrho_0,u_0,\theta_0)\|_{H^2}^2   \label{higneuestimate}
\end{eqnarray}
for any $t\geq 0$,
where $\Lambda$ is constructed by \eqref{growth}, and the constant $C$ may depend on $g$, $\mu$,
$\mu_0$, $\bar{e}$, $\bar{\rho}$, $\Lambda$ and $\Omega$.
\end{pro}
\begin{pf}
The first estimate \eqref{uestimate} can be shown by an argument similar to that in Lemma \ref{erroestimate}.
In fact, following the process in the derivation of \eqref{12safd} and \eqref{inequalemee}, we obtain the following two inequalities
 \begin{equation}\label{Ainequfirs}
\begin{aligned}
& \frac{\mm{d}}{\mm{d}t}\|\sqrt{\bar{\rho}} u(t)\|^2_{L^2}+  {\mu}\|\nabla u(t)\|^2_{L^2}+
 {\mu_0}\|\mm{div} u(t)\|^2_{L^2} \\
&\leq I_1+2\Lambda\left[ \|\sqrt{\bar{\rho}} u\|^2_{L^2} +\int_0^t(\mu\|\nabla u\|^2_{L^2}
+ {\mu}_0\|\mm{div} u\|^2_{L^2}) \mathrm{d}\tau\right]
\end{aligned}
\end{equation}
and \begin{equation}\label{inequaldf2}
\begin{aligned}
&\frac{1}{\Lambda}\|\sqrt{\bar{\rho}}  u_t(t)\|^2_{L^2}+ {\mu}\|\nabla u(t)\|^2_{L^2}+
 {\mu_0}\|\mm{div} u(t)\|^2_{L^2} \\
& \leq  I_1+ {\Lambda}\|\sqrt{\bar{\rho}} u(t)\|^2_{L^2}+2{\Lambda}\int_0^t( {\mu}\|\nabla u\|^2_{L^2}+
 {\mu_0}\|\mm{div} u\|^2_{L^2})\mm{d}\tau
\end{aligned}\end{equation}
with $$I_1=\left\{2(\mu\|\nabla u\|^2_{L^2} +{\mu_0}\|\mm{div} u\|_{L^2}^2)+\frac{1}{\Lambda}\int \left\{\bar{\rho}| u_t|^2
-g\bar{\rho}'{u}_3^2+
[(1+a)\bar{p}\mm{div} u
-2g\bar{\rho}{u}_3]\mm{div}{ u}\right\}\mathrm{d} x\right\}\bigg|_{t=0}.$$
 An application of Gronwall's inequality to \eqref{Ainequfirs} implies that for any $t\geq 0$,
 \begin{equation*} \begin{aligned}\label{}
\|\sqrt{\bar{\rho}} u(t)\|^2_{L^2}+ \int_0^t(\mu\|\nabla u\|^2_{L^2}+{\mu}_0\|\mm{div} u\|^2_{L^2})\mm{d}\tau
\leq & e^{2\Lambda t}\|\sqrt{\bar{\rho}} u_0\|^2_{L^2} +\frac{I_1}{2\Lambda}\left(e^{2\Lambda t}-1\right) \\
\leq & Ce^{2\Lambda t}(\|(\varrho_0,\theta_0)\|_{H^1}^2+ \| u_0\|_{H^2}^2),
\end{aligned}\end{equation*}
which, together with Poincar\'e's inequality and \eqref{inequaldf2}, results in
\begin{eqnarray*}\label{}
\| u(t)\|_{H^1 }^2+\| u_t(t)\|^2_{L^2 }+ \int_0^t\|\nabla u(s)\|^2_{L^2}\mm{d}s
\leq Ce^{2\Lambda t}(\|(\varrho_0,\theta_0)\|_{H^1}^2+\| u_0\|_{H^2}^2).
\end{eqnarray*}
Thus, using \eqref{0108}$_1$ and \eqref{0108}$_2$, we have
\begin{equation*}\begin{aligned}  \label{}
\|(\varrho,\theta)(t)\|_{L^2}\leq & \|(\varrho_0,\theta_0)\|_{L^2}+\int_0^t
\|(\varrho,\theta)_s(s)\|_{L^2}\mm{d}s \\
\leq &\|(\varrho_0,\theta_0)\|_{L^2}+(1+a)\|(\bar{\rho},\bar{\theta})
\|_{H^1}\int_0^t\| u(s)\|_{H^1}\mm{d}s \\
\leq & Ce^{\Lambda t}(\|(\varrho_0,\theta_0)\|_{H^1}+\| u_0\|_{H^2}).
\end{aligned}\end{equation*}
 Hence the estimate \eqref{uestimate} follows from the above two estimates.

Finally, following the arguments in the proof of \eqref{ernygfor}, one find that
\begin{equation*}\label{}\begin{aligned}
& \|(\varrho,u,\theta)(t)\|_{H^2}^2+\int_0^t\left[\|
 u\|_{H^3}^2+\left\| (\bar{e}\rho_t+\bar{\rho}\theta_t)\right\|_{H^2}^2\right]\mm{d}\tau  \\
& \leq C\|(\varrho_0,u_0,\theta_0)\|_{H^2}^2+ \int_0^t[ C\|(\varrho,u,\theta)\|_{L^2}^2
+\Lambda\|(\varrho,u,\theta)\|_{H^2}^2]\mm{d}\tau,
\end{aligned}
\end{equation*}
 which, combined with \eqref{uestimate}, gives \eqref{higneuestimate} due to Gronwall's inequality. This completes the proof.
  \hfill $\Box$
\end{pf}

\vspace{4mm} \noindent\textbf{Acknowledgements.}
The research of Fei Jiang was supported by NSFC (Grant No. 11671086), the NSF of Fujian Province of China (Grant No. 2016J06001) and the Education Department of Fujian Province (Grant No. SX2015-02).
 Finally, the authors would like to thank the anonymous referee for invaluable suggestions,
  which improve the presentation of this paper.

\renewcommand\refname{References}
\renewenvironment{thebibliography}[1]{%
\section*{\refname}
\list{{\arabic{enumi}}}{\def\makelabel##1{\hss{##1}}\topsep=0mm
\parsep=0mm
\partopsep=0mm\itemsep=0mm
\labelsep=1ex\itemindent=0mm
\settowidth\labelwidth{\small[#1]}%
\leftmargin\labelwidth \advance\leftmargin\labelsep
\advance\leftmargin -\itemindent
\usecounter{enumi}}\small
\def\newblock{\ }
\sloppy\clubpenalty4000\widowpenalty4000
\sfcode`\.=1000\relax}{\endlist}
\bibliographystyle{model1b-num-names}

\begin{thebibliography}{47}
\expandafter\ifx\csname natexlab\endcsname\relax\def\natexlab#1{#1}\fi
\providecommand{\bibinfo}[2]{#2}
\ifx\xfnm\relax \def\xfnm[#1]{\unskip,\space#1}\fi
\bibitem[{Adams and John(2005)}]{ARAJJFF}
\bibinfo{author}{R.A. Adams}, \bibinfo{author}{J.F.F. John},
  \bibinfo{title}{{Sobolev Space}}, \bibinfo{publisher}{Academic Press: New
  York}, \bibinfo{year}{2005}.
\bibitem[{Aye and Nishida£¬T.(1998)}]{APTNHCC}
\bibinfo{author}{P.~Aye}, \bibinfo{author}{Nishida£¬T.}, \bibinfo{title}{{Heat
  convection of compressible fluid. In: Fujita, H.(ed) Recent Developments in
  Domain Decomposition Methods and Flow Problems}},
  \bibinfo{journal}{Mathematical Sciences and Applications, Gakkotosho, Tokyo,
  Japan} \bibinfo{volume}{11} (\bibinfo{year}{1998}) \bibinfo{pages}{107--115}.
\bibitem[{Benardand(1900)}]{BHERESS}
\bibinfo{author}{H.~Benardand}, \bibinfo{title}{Les tourbillons cellulaires
  dans une nappe liquide}, \bibinfo{journal}{Revue g\'en\'erale des Sciences
  pures et appliqu\'ees} \bibinfo{volume}{45} (\bibinfo{year}{1900})
  \bibinfo{pages}{1261--71 and 1309--28}.
\bibitem[{Cattaneo et~al.(1991)Cattaneo, Brummell, Toomre, Malagoli and
  Hurlburt}]{CFBNHTJMAHAE}
\bibinfo{author}{F.~Cattaneo}, \bibinfo{author}{N.~Brummell},
  \bibinfo{author}{J.~Toomre}, \bibinfo{author}{A.~Malagoli},
  \bibinfo{author}{N.~Hurlburt}, \bibinfo{title}{Turbulent compressible
  convection}, \bibinfo{journal}{The Astrophysical Journal}
  \bibinfo{volume}{370} (\bibinfo{year}{1991}) \bibinfo{pages}{282--294}.
\bibitem[{Chandrasekhar(1939)}]{CSAiTTS}
\bibinfo{author}{S.~Chandrasekhar}, \bibinfo{title}{{An Introduction to the
  Study of Stellar Structures}}, \bibinfo{publisher}{University of Chicago
  Press}, \bibinfo{year}{1939}.
\bibitem[{Chandrasekhar(1961)}]{CSHHSCPO}
\bibinfo{author}{S.~Chandrasekhar}, \bibinfo{title}{{Hydrodynamic and
  Hydromagnetic Stability, The International Series of Monographs on Physics}},
  \bibinfo{publisher}{Oxford, Clarendon Press}, \bibinfo{year}{1961}.
\bibitem[{Coscia and Padula(1990)}]{CosciaPadula}
\bibinfo{author}{V.~Coscia}, \bibinfo{author}{M.~Padula},
  \bibinfo{title}{{Nonlinear energy stability in a compressible atmosphere}},
  \bibinfo{journal}{Geophys. Astrophys. Fluid. Dynmics} \bibinfo{volume}{54}
  (\bibinfo{year}{1990}) \bibinfo{pages}{49--83}.
\bibitem[{Drazin and Reid(2004)}]{DPGRWHHC}
\bibinfo{author}{P.G. Drazin}, \bibinfo{author}{W.H. Reid},
  \bibinfo{title}{{Hydrodynamic Stability, 2nd}}, \bibinfo{publisher}{Cambridge
  University Press}, \bibinfo{year}{2004}.
\bibitem[{Field(1965)}]{FGBFE}
\bibinfo{author}{G.B. Field}, \bibinfo{title}{{ Thermal instability}},
  \bibinfo{journal}{Astrophysical Journal} \bibinfo{volume}{142}
  (\bibinfo{year}{1965}) \bibinfo{pages}{531--567}.
\bibitem[{Friedlander et~al.(1997)Friedlander, Strauss and Vishik}]{FSSWVMNA}
\bibinfo{author}{S.~Friedlander}, \bibinfo{author}{W.~Strauss},
  \bibinfo{author}{M.~Vishik}, \bibinfo{title}{{Nonlinear instability in an
  ideal fluid}}, \bibinfo{journal}{Annales de l'Institut Henri Poincare (C) Non
  Linear Analysis} \bibinfo{volume}{14} (\bibinfo{year}{1997})
  \bibinfo{pages}{187--209}.
\bibitem[{Friedlander and Vishik(2003)}]{FrrVishikM}
\bibinfo{author}{S.~Friedlander}, \bibinfo{author}{M.~Vishik},
  \bibinfo{title}{{Nonlinear instability in two dimensional ideal fluids: the
  case of a dominant eigenvalue}}, \bibinfo{journal}{Comm. Math. Phys.}
  \bibinfo{volume}{243} (\bibinfo{year}{2003}) \bibinfo{pages}{261--273}.
\bibitem[{Galdi(1984)}]{GGPTA}
\bibinfo{author}{G.P. Galdi}, \bibinfo{title}{{ The rotating B\'enard problem:
  a nonlinear energy stability analysis. (Rome, 1984), Teubner, Stuttgart},},
  \bibinfo{journal}{Applications of mathematics in technology}
  (\bibinfo{year}{1984}) \bibinfo{pages}{79--95}.
\bibitem[{Galdi(1985)}]{GGPNA}
\bibinfo{author}{G.P. Galdi}, \bibinfo{title}{{Nonlinear stability of the
  magnetic B\'enard problem via a generalized energy method}},
  \bibinfo{journal}{Arch. Rational Mech. Anal.} \bibinfo{volume}{62(2)}
  (\bibinfo{year}{1985}) \bibinfo{pages}{167--186}.
\bibitem[{Galdi and Padula(1989)}]{GGPPMNA}
\bibinfo{author}{G.P. Galdi}, \bibinfo{author}{M.~Padula}, \bibinfo{title}{{New
  contributions to nonlinear stability of the magnetic B\'enard problem}},
  \bibinfo{journal}{Applications of mathematics in industry and technology
  (Siena, 1988), Teubner, Stuttgart}  (\bibinfo{year}{1989})
  \bibinfo{pages}{166--178}.
\bibitem[{Galdi and Padula(1991)}]{GGPPMFR}
\bibinfo{author}{G.P. Galdi}, \bibinfo{author}{M.~Padula},
  \bibinfo{title}{{Further results in the nonlinear stability of the magnetic
  B\'enard problem}}, \bibinfo{journal}{Mathematical aspects of fluid and
  plasma dynamics (Salice Terme, 1988), Lecture Notes in Math., 1460, Springer,
  Berlin}  (\bibinfo{year}{1991}).
\bibitem[{Galdi and Straughan(1985)}]{GGPSBAP}
\bibinfo{author}{G.P. Galdi}, \bibinfo{author}{B.~Straughan},
  \bibinfo{title}{{A nonlinear analysis of the stabilizing effect of rotation
  in the B\'enard problem }}, \bibinfo{journal}{Proc. Roy. Soc. London Ser. A}
  \bibinfo{volume}{402(1823)} (\bibinfo{year}{1985}) \bibinfo{pages}{257--283}.
\bibitem[{Gough et~al.(1976)Gough, Moore, Spiegel and Weiss}]{GDOMDRSEAWNO}
\bibinfo{author}{D.O. Gough}, \bibinfo{author}{D.R. Moore},
  \bibinfo{author}{E.A. Spiegel}, \bibinfo{author}{N.O. Weiss},
  \bibinfo{title}{{ Convective instability in a compressible atmosphere II}},
  \bibinfo{journal}{Ap. J.} \bibinfo{volume}{206} (\bibinfo{year}{1976})
  \bibinfo{pages}{536--542}.
\bibitem[{Graham(1975)}]{GENJFM}
\bibinfo{author}{E.~Graham}, \bibinfo{title}{{Numerical simulation of
  two-dimensional compressible convection}}, \bibinfo{journal}{J. Fluid Mech.}
  \bibinfo{volume}{70} (\bibinfo{year}{1975}) \bibinfo{pages}{689}.
\bibitem[{Guidoboni and Padula(2005)}]{GGMPO}
\bibinfo{author}{G.~Guidoboni}, \bibinfo{author}{M.~Padula},
  \bibinfo{title}{{On the B\'enard problem}}, \bibinfo{journal}{Progress in
  Nonlinear Differential Equations and Their Applications} \bibinfo{volume}{61}
  (\bibinfo{year}{2005}) \bibinfo{pages}{137--148}.
\bibitem[{Guo et~al.(2007)Guo, Hallstrom and Spirn}]{GYHCSDDC}
\bibinfo{author}{Y.~Guo}, \bibinfo{author}{C.~Hallstrom},
  \bibinfo{author}{D.~Spirn}, \bibinfo{title}{{Dynamics near unstable,
  interfacial fluids}}, \bibinfo{journal}{Commun. Math. Phys.}
  \bibinfo{volume}{270} (\bibinfo{year}{2007}) \bibinfo{pages}{635--689}.
\bibitem[{Guo and Han(2010)}]{GYYQH}
\bibinfo{author}{Y.~Guo}, \bibinfo{author}{Y.Q. Han}, \bibinfo{title}{{Critical
  Rayleigh number in Rayleigh-B$\mm{\acute{e}}$nard convection}},
  \bibinfo{journal}{Quart. Appl. Math.} \bibinfo{volume}{LXVIII}
  (\bibinfo{year}{2010}) \bibinfo{pages}{149--160}.
\bibitem[{Guo and Strauss(1995)}]{GYSWIC}
\bibinfo{author}{Y.~Guo}, \bibinfo{author}{W.~Strauss},
  \bibinfo{title}{{Instability of periodic BGK equilibria}},
  \bibinfo{journal}{Comm. Pure Appl. Math.} \bibinfo{volume}{48}
  (\bibinfo{year}{1995}) \bibinfo{pages}{861--894}.
\bibitem[{Guo and Tice(2011)}]{GYTI2}
\bibinfo{author}{Y.~Guo}, \bibinfo{author}{I.~Tice}, \bibinfo{title}{Linear
  rayleigh-taylor instability for viscous, compressible fluids},
  \bibinfo{journal}{SIAM J. Math. Anal.} \bibinfo{volume}{42}
  (\bibinfo{year}{2011}) \bibinfo{pages}{1688--1720}.
\bibitem[{Hwang and Guo(2003)}]{HHJGY}
\bibinfo{author}{H.J. Hwang}, \bibinfo{author}{Y.~Guo}, \bibinfo{title}{{On the
  dynamical Rayleigh-Taylor instability}}, \bibinfo{journal}{Arch. Rational
  Mech. Anal.} \bibinfo{volume}{167} (\bibinfo{year}{2003})
  \bibinfo{pages}{235--253}.
\bibitem[{Jang and Tice(2011)}]{JJTIIA}
\bibinfo{author}{J.~Jang}, \bibinfo{author}{I.~Tice},
  \bibinfo{title}{{Instability theory of the Navier-Stokes-Poisson equations}},
  \bibinfo{journal}{To appear in Analysis $\&$ PDE}  (\bibinfo{year}{2011}).
\bibitem[{Jeffreys(1930)}]{JHTICFHB}
\bibinfo{author}{H.~Jeffreys}, \bibinfo{title}{{The instability of a
  compressible fluid heated below}}, \bibinfo{journal}{Proc. Cambridge Phil.
  Soc.} \bibinfo{volume}{26} (\bibinfo{year}{1930}) \bibinfo{pages}{170--172}.
\bibitem[{Jiang and Jiang(2014)}]{JFJSO2014}
\bibinfo{author}{F.~Jiang}, \bibinfo{author}{S.~Jiang}, \bibinfo{title}{{On
  instability and stability of three-dimensional gravity flows in a bounded
  domain}}, \bibinfo{journal}{Adv. Math.} \bibinfo{volume}{264}
  (\bibinfo{year}{2014}) \bibinfo{pages}{831--863}.
\bibitem[{Jiang et~al.(2014)Jiang, Jiang and Wang}]{JFJSWWWOA}
\bibinfo{author}{F.~Jiang}, \bibinfo{author}{S.~Jiang}, \bibinfo{author}{Y.J.
  Wang}, \bibinfo{title}{{On the Rayleigh-Taylor instability for the
  incompressible viscous magnetohydrodynamic equations}},
  \bibinfo{journal}{Comm. Partial Differential Equations} \bibinfo{volume}{39}
  (\bibinfo{year}{2014}) \bibinfo{pages}{399--438}.
\bibitem[{Joseph(1966)}]{JDDN1966}
\bibinfo{author}{D.D. Joseph}, \bibinfo{title}{{Nonlinear stability of the
  Boussinesq equations by the method of energy}}, \bibinfo{journal}{Arch.
  Rational Mech. Anal.} \bibinfo{volume}{22} (\bibinfo{year}{1966})
  \bibinfo{pages}{163--184}.
\bibitem[{Kaniel(1967)}]{KSKA201710110936}
\bibinfo{author}{S. Kaniel},
\bibinfo{author}{A. Kovetz},
 \bibinfo{title}{{Schwarzschild's criterion for instability}}, \bibinfo{journal}{The Physics of Fluids} \bibinfo{volume}{10} (\bibinfo{year}{1967})
  \bibinfo{pages}{1186--1193}.
\bibitem[{Kawashima(1983)}]{kawashima1984systems}
\bibinfo{author}{S.~Kawashima}, \bibinfo{title}{Systems of a
  hyperbolic-parabolic composite type, with applications to the equations of
  magnetohydrodynamics}, \bibinfo{publisher}{Ph. D. Thesis, Kyoto University},
  \bibinfo{year}{1983}.
\bibitem[{Lebovitz(1986)}]{LEbovitzNROth}
\bibinfo{author}{N.R. Lebovitz}, \bibinfo{title}{On the necessity of
  schwarzshild's criterion for stability}, \bibinfo{journal}{Astrophysical
  Journal} \bibinfo{volume}{146} (\bibinfo{year}{1986})
  \bibinfo{pages}{946--949}.
\bibitem[{Ma and Wang(2004)}]{MTWSB}
\bibinfo{author}{T.~Ma}, \bibinfo{author}{S.~Wang}, \bibinfo{title}{Dynamic
  bifurcation and stability in the rayleigh--b\'enard convection},
  \bibinfo{journal}{Comm. Math. Sci} \bibinfo{volume}{2} (\bibinfo{year}{2004})
  \bibinfo{pages}{159--183}.
\bibitem[{Matsumura and Nishida(1979)}]{MANTTP351}
\bibinfo{author}{A.~Matsumura}, \bibinfo{author}{T.~Nishida},
  \bibinfo{title}{The initial value problem for the equation of compressible
  viscous and heat-conductive fluids}, \bibinfo{journal}{Proc. Jpn. Acad.
  Ser-A.} \bibinfo{volume}{55} (\bibinfo{year}{1979})
  \bibinfo{pages}{337--342}.
\bibitem[{Matsumura and Nishida(1982)}]{MANTIJ321}
\bibinfo{author}{A.~Matsumura}, \bibinfo{author}{T.~Nishida},
  \bibinfo{title}{Initial-boundary value problems for the equations of motion
  of general fluids, computing methods in applied sciences and engineering},
  \bibinfo{journal}{J. Math. Kyoto. Univ. V (Versailles, 1981)}
  \bibinfo{volume}{20} (\bibinfo{year}{1982}) \bibinfo{pages}{389--406}.
\bibitem[{Matsumura and Nishida(1983)}]{MANTIC481}
\bibinfo{author}{A.~Matsumura}, \bibinfo{author}{T.~Nishida},
  \bibinfo{title}{Initial boundary value problems for the equations of motion
  of compressible viscous and heat conductive fluids}, \bibinfo{journal}{Comm.
  Math. Phys.} \bibinfo{volume}{89} (\bibinfo{year}{1983})
  \bibinfo{pages}{445--464}.
\bibitem[{Nishida et~al.(2012)Nishida, Padula and Teramoto}]{TNMPYTHCI}
\bibinfo{author}{T.~Nishida}, \bibinfo{author}{M.~Padula},
  \bibinfo{author}{Y.~Teramoto}, \bibinfo{title}{{ Heat Convection of
  Compressible Viscous Fluids. I}}, \bibinfo{journal}{J. Math. Fluid Mech.}
  \bibinfo{volume}{15} (\bibinfo{year}{2012}) \bibinfo{pages}{525--536}.
\bibitem[{Nishida et~al.(2013)Nishida, Padula and Teramoto}]{TNMPYTHCII}
\bibinfo{author}{T.~Nishida}, \bibinfo{author}{M.~Padula},
  \bibinfo{author}{Y.~Teramoto}, \bibinfo{title}{{ Heat Convection of
  Compressible Viscous Fluids. II}}, \bibinfo{journal}{J. Math. Fluid Mech.}
  \bibinfo{volume}{15} (\bibinfo{year}{2013}) \bibinfo{pages}{689--700}.
\bibitem[{Novotn{\`y} and Stra{\v{s}}kraba(2004)}]{NASII04}
\bibinfo{author}{A.~Novotn{\`y}}, \bibinfo{author}{I.~Stra{\v{s}}kraba},
  \bibinfo{title}{{Introduction to the Mathematical Theory of Compressible
  Flow}}, \bibinfo{publisher}{Oxford University Press, USA},
  \bibinfo{year}{2004}.
\bibitem[{Padula and Bollettmo(1986)}]{OMBUMINst}
\bibinfo{author}{M.~Padula}, \bibinfo{author}{U.M.I. Bollettmo},
  \bibinfo{title}{{ Nonlinear energy stability for the compressible B\'enard
  problem}}, \bibinfo{journal}{Boll. Un. Mat. Ital. B} \bibinfo{volume}{5B}
  (\bibinfo{year}{1986}) \bibinfo{pages}{581--602}.
\bibitem[{Rayleigh(1916)}]{RALOC}
\bibinfo{author}{L.~Rayleigh}, \bibinfo{title}{{On convective currants in a
  horizontal layer of fluid when the higher termperature is on the under
  side}}, \bibinfo{journal}{Phil. Mag.} \bibinfo{volume}{32}
  (\bibinfo{year}{1916}) \bibinfo{pages}{529--546}.
\bibitem[{Rosencrans(1969)}]{SIAMJAMOR}
\bibinfo{author}{S.~Rosencrans}, \bibinfo{title}{On schwarzschild's criterion},
  \bibinfo{journal}{SIAM J. Appl. Math.} \bibinfo{volume}{17}
  (\bibinfo{year}{1969}) \bibinfo{pages}{231--239}.
\bibitem[{Rumford(1870)}]{CROFTCW}
\bibinfo{author}{C.~Rumford}, \bibinfo{title}{{Of the propagation of heat in
  fluids, Complete Works, 1, 239}}, \bibinfo{publisher}{American Academy of
  Arts and Sciences, Boston}, \bibinfo{year}{1870}.
\bibitem[{Schwarzschild(1906)}]{KSNKGWG}
\bibinfo{author}{K.~Schwarzschild}, \bibinfo{title}{\"uber das gleichgewicht
  der sonnenatmosph\"are}, \bibinfo{journal}{Nachr. Kgl. Ges. Wiss. G\"ottingen
  Math.-Phys. Klasse}  (\bibinfo{year}{1906}) \bibinfo{pages}{41--53}.
\bibitem[{Schwarzschild(1961)}]{MSCTAJ}
\bibinfo{author}{M.~Schwarzschild}, \bibinfo{title}{{Convection in stars}},
  \bibinfo{journal}{The Astrophysical Journal} \bibinfo{volume}{134(1)}
  (\bibinfo{year}{1961}) \bibinfo{pages}{1--8}.
\bibitem[{Spiegel(1965)}]{SEAAPKC}
\bibinfo{author}{E.A. Spiegel}, \bibinfo{title}{Convective instability in a
  compressible atmosphere i}, \bibinfo{journal}{Ap. J.} \bibinfo{volume}{141}
  (\bibinfo{year}{1965}) \bibinfo{pages}{1068--1090}.
\bibitem[{Thompson(1882)}]{TJOACteins}
\bibinfo{author}{J.~Thompson}, \bibinfo{title}{{On a changing
  tesselated.structure in certain liquids}}, \bibinfo{journal}{Pro. Phil. Soc.
  Glasgow} \bibinfo{volume}{13} (\bibinfo{year}{1882})
  \bibinfo{pages}{464--468}.
\bibitem[{Unno et~al.(1960)Unno, Kato and Makita}]{UWKSMMC}
\bibinfo{author}{W.~Unno}, \bibinfo{author}{S.~Kato},
  \bibinfo{author}{M.~Makita}, \bibinfo{title}{{Convective instability in
  polytropic atmosphere. I}}, \bibinfo{journal}{Pub. Astr. Soc. Japan}
  \bibinfo{volume}{12(2)} (\bibinfo{year}{1960}) \bibinfo{pages}{192--202}.

\end{thebibliography}

\end{document}